
\input amssym.def
\input amssym.tex


\def\item#1{\vskip1.3pt\hang\textindent {\rm #1}}

\def\itemitem#1{\vskip1.3pt\indent\hangindent2\parindent\textindent {\rm #1}}

\tolerance=300
\pretolerance=200
\hfuzz=1pt
\vfuzz=1pt


\hoffset=0.6in
\voffset=0.8in

\hsize=5.8 true in 


\vsize=8.5 true in
\parindent=25pt
\mathsurround=1pt
\parskip=1pt plus .25pt minus .25pt
\normallineskiplimit=.99pt

\countdef\revised=100
\mathchardef\emptyset="001F 
\chardef\ss="19
\def\3{\ss}
\def\anf{$\lower1.2ex\hbox{"}$}
\def\frac#1#2{{#1 \over #2}}
\def\>{>\!\!>}
\def\<{<\!\!<}

\def\into{\hookrightarrow}
\def\ssarr{\hbox to 30pt{\rightarrowfill}}
\def\sarr{\hbox to 40pt{\rightarrowfill}}
\def\arr{\hbox to 60pt{\rightarrowfill}}
\def\larr{\hbox to 60pt{\leftarrowfill}}
\def\Arr{\hbox to 80pt{\rightarrowfill}}

{}

\def\ad{\mathop{\rm ad}\nolimits}

\def\Ad{\mathop{\rm Ad}\nolimits}

\def\Aut{\mathop{\rm Aut}\nolimits}

\def\conv{\mathop{\rm conv}\nolimits}

\def\det{\mathop{\rm det}\nolimits}
\def\diag{\mathop{\rm diag}\nolimits}

\def\End{\mathop{\rm End}\nolimits}
\def\Ext{\mathop{\rm Ext}\nolimits}

\def\Gl{\mathop{\rm Gl}\nolimits}
\def\GL{\mathop{\rm GL}\nolimits}

\def\Herm{\mathop{\rm Herm}\nolimits}

%
%

\def\Im{\mathop{\rm Im}\nolimits}

\def\Int{\mathop{\rm int}\nolimits}


\def\sgn{\mathop{\rm sgn}\nolimits}
\def\Sl{\mathop{\rm Sl}\nolimits}
\def\SO{\mathop{\rm SO}\nolimits}

\def\Symm{\mathop{\rm Symm}\nolimits}
\def\Sp{\mathop{\rm Sp}\nolimits}
\def\Spec{\mathop{\rm Spec}\nolimits}

\def\SU{\mathop{\rm SU}\nolimits}
\def\sup{\mathop{\rm sup}\nolimits}



\def\0{{\bf 0}}
\def\1{{\bf 1}}

\def\a{{\frak a}}

\def\e{{\frak e}}
\def\f{{\frak f}}
\def\g{{\frak g}}
\def\gl{{\frak {gl}}}
\def\h{{\frak h}}

\def\k{{\frak k}}
\def\l{{\frak l}}
\def\m{{\frak m}}

\def\n{{\frak n}}

\def\p{{\frak p}}
\def\q{{\frak q}}

\def\s{{\frak s}}

\def\sp{{\frak {sp}}}

\def\su{{\frak {su}}}
\def\so{{\frak {so}}}
\def\sL{{\frak {sl}}}

\def\uu{{\frak u}}

\def\z{{\frak z}}

\def\C{{\Bbb C}} 
\def\D{{\Bbb D}} 
 
\def\H{{\Bbb H}}

\def\N{{\Bbb N}}

\def\R{{\Bbb R}} 
 
\def\Z{{\Bbb Z}} 

\def\:{\colon}  
\def\.{{\cdot}}
\def\|{\Vert}
\def\bsk{\bigskip}

\def\giantskip{\vskip2\bigskipamount}
\def\gsk{\giantskip}
\def \la {\langle}
\def\msk{\medskip}
\def \ra {\rangle}
\def \res {\!\mid\!\!}

\def\ssk{\smallskip}

\def\bbr{\bigbreak}
\def\giantbreak{\par \ifdim\lastskip<2\bigskipamount \removelastskip
         \penalty-400 \giantskip\fi}

\def\nin{\noindent}
\def\cen{\centerline}
\def\pagebreak{\vskip 0pt plus 0.0001fil\break}
\def\linebreak{\break}

\def\hat{\widehat}

\def\eps{\varepsilon}
\def\phi{\varphi}
\def\epsilon{\varepsilon}

\def\nin{\noindent}
\def\oline{\overline}

\def\pder#1,#2,#3 { {\partial #1 \over \partial #2}(#3)}
\def\pde#1,#2 { {\partial #1 \over \partial #2}}


\def\subeq{\subseteq}
\def\supeq{\supseteq}

\def\up{{\uparrow}}

\font\eightrm=cmr8


\font\smc=cmcsc10
\font\bfone=cmbx10 scaled\magstep1 
\font\bftwo=cmbx10 scaled\magstep2 

\def\qed{{\unskip\nobreak\hfil\penalty50\hskip .001pt \hbox{}\nobreak\hfil
          \vrule height 1.2ex width 1.1ex depth -.1ex
           \parfillskip=0pt\finalhyphendemerits=0\medbreak}\rm}

\def\qeddis{\eqno{\vrule height 1.2ex width 1.1ex depth -.1ex} $$
                   \medbreak\rm}

\def\Lemma #1. {\bigbreak\vskip-\parskip\noindent{\bf Lemma #1.}\quad\it}

\def\Sublemma #1. {\bigbreak\vskip-\parskip\noindent{\bf Sublemma #1.}\quad\it}

\def\Proposition #1. {\bigbreak\vskip-\parskip\noindent{\bf Proposition #1.}
\quad\it}

\def\Corollary #1. {\bigbreak\vskip-\parskip\nin{\bf Corollary #1.}
\quad\it}

\def\Theorem #1. {\bigbreak\vskip-\parskip\noindent{\bf Theorem #1.}
\quad\it}

\def\Definition #1. {\rm\bigbreak\vskip-\parskip\noindent{\bf Definition #1.}
\quad}

\def\Remark #1. {\rm\bigbreak\vskip-\parskip\noindent{\bf Remark #1.}\quad}

\def\Example #1. {\rm\bigbreak\vskip-\parskip\noindent{\bf Example #1.}\quad}

\def\Problems #1. {\bigbreak\vskip-\parskip\noindent{\bf Problems #1.}\quad}
\def\Problem #1. {\bigbreak\vskip-\parskip\noindent{\bf Problems #1.}\quad}

\def\Conjecture #1. {\bigbreak\vskip-\parskip\noindent{\bf Conjecture #1.}\quad}

\def\Proof#1.{\rm\par\ifdim\lastskip<\bigskipamount\removelastskip\fi\smallskip
            \noindent {\bf Proof.}\quad}

\def\Axiom #1. {\bigbreak\vskip-\parskip\noindent{\bf Axiom #1.}\quad\it}

\def\Satz #1. {\bigbreak\vskip-\parskip\noindent{\bf Satz #1.}\quad\it}

\def\Korollar #1. {\bbr\vskip-\parskip\nin{\bf Korollar #1.} \quad\it}

\def\Bemerkung #1. {\rm\bigbreak\vskip-\parskip\noindent{\bf Bemerkung #1.}
\quad}

\def\Beispiel #1. {\rm\bigbreak\vskip-\parskip\noindent{\bf Beispiel #1.}\quad}
\def\Aufgabe #1. {\rm\bigbreak\vskip-\parskip\noindent{\bf Aufgabe #1.}\quad}

\def\Beweis#1. {\rm\par\ifdim\lastskip<\bigskipamount\removelastskip\fi
           \smallskip\noindent {\bf Beweis.}\quad}

\nopagenumbers

\def\date{\ifcase\month\or January\or February \or March\or April\or May
\or June\or July\or August\or September\or October\or November
\or December\fi\space\number\day, \number\year}

\def\title{Title ??}
\def\author{Author ??}

\def\thanks#1{\footnote*{\eightrm#1}}

\def\rightheadline{\hfil{\eightrm\title}\hfil\tenbf\folio}
\def\leftheadline{\tenbf\folio\hfil{\eightrm\author}\hfil}
\headline={\vbox{\line{\ifodd\pageno\rightheadline\else\leftheadline\fi}}}

\def\firstheadline{}
\def\firstfootline{\cen{\rm\folio}}

\def\seite #1 {\pageno #1
               \headline={\ifnum\pageno=#1 \firstheadline
               \else\ifodd\pageno\rightheadline\else\leftheadline\fi\fi}
               \footline={\ifnum\pageno=#1 \firstfootline\else{}\fi}}

\newdimen\dimenone
 \def\checkleftspace#1#2#3#4{
 \dimenone=\pagetotal
 \advance\dimenone by -\pageshrink   
 \ifdim\dimenone>\pagegoal          
   \else\dimenone=\pagetotal
        \advance\dimenone by \pagestretch
        \ifdim\dimenone<\pagegoal
          \dimenone=\pagetotal
          \advance\dimenone by#1         
          \setbox0=\vbox{#2\parskip=0pt                
                     \hyphenpenalty=10000
                     \rightskip=0pt plus 5em
                     \noindent#3 \vskip#4}    
        \advance\dimenone by\ht0
        \advance\dimenone by 3\baselineskip   
        \ifdim\dimenone>\pagegoal\vfill\eject\fi
          \else\eject\fi\fi}


\def\subheadline #1{\nin\bigbreak\vskip-\lastskip
      \checkleftspace{0.7cm}{\bf}{#1}{\medskipamount}
          \indent\vskip0.7cm\centerline{\bf #1}\medskip}

\def\sectionheadline #1{\bigbreak\vskip-\lastskip
      \checkleftspace{1.1cm}{\bf}{#1}{\bigskipamount}
         \vbox{\vskip1.1cm}\cen{\bfone #1}\bsk}

\def\lsectionheadline #1 #2{\bigbreak\vskip-\lastskip
      \checkleftspace{1.1cm}{\bf}{#1}{\bigskipamount}
         \vbox{\vskip1.1cm}\cen{\bfone #1}\msk \cen{\bfone #2}\bsk}

\def\lchapterheadline #1 #2{\bigbreak\vskip-\lastskip\indent\vskip3cm
                       \cen{\bftwo #1} \msk \cen{\bftwo #2} \gsk}
\def\llsectionheadline #1 #2 #3{\bigbreak\vskip-\lastskip\indent\vskip1.8cm
\cen{\bfone #1} \msk \cen{\bfone #2} \msk \cen{\bfone #3} \nobreak\bsk\nobreak}


\newtoks\literat
\def\[#1 #2\par{\literat={#2\unskip.}%
\hbox{\vtop{\hsize=.15\hsize\nin [#1]\hfill}
\vtop{\hsize=.82\hsize\nin\the\literat}}\par
\vskip.3\baselineskip}

\mathchardef\emptyset="001F 
\def\address{Author: \tt$\backslash$def$\backslash$address$\{$??$\}$}

\def\firstpage{\nin
{\obeylines \parindent 0pt }
\vskip2cm
\centerline {\bfone \title}
\gsk
\centerline{\bf\author}

\vskip1.5cm \rm}

\def\addresstwo{}

\def\dlastpage{\par\vbox{\vskip1cm\nin
\line{
\vtop{\hsize=.5\hsize{\parindent=0pt\baselineskip=10pt\nin\address}}
\quad 
\vtop{\hsize=.42\hsize\nin{\parindent=0pt
\baselineskip=10pt\addresstwo}}
\hfill} }}


\def\firstpage{\nin
{\obeylines \parindent 0pt }
\vskip2cm
\centerline {\bfone \title}
\ssk
\centerline {\bfone \titletwo}
\gsk
\centerline{\bf\author}
\vskip1.5cm \rm}
\def\bs{\backslash} 
\def\addots{\mathinner{\mkern1mu\raise1pt\vbox{\kern7pt\hbox{.}}\mkern2mu
\raise4pt\hbox{.}\mkern2mu\raise7pt\hbox{.}\mkern1mu}}

\pageno=1
\def\up#1{\leavevmode \raise.16ex\hbox{#1}}
 at 8truept
 at 8truept
 at 12truept
\chardef\ss="19
\def\3{\ss}
\def\title{Complex crowns of Riemannian symmetric spaces}
\def\titletwo{and non-compactly causal symmetric spaces}
\def\author{Simon Gindikin$^*$ and Bernhard Kr\"otz${}^\dagger$}
\footnote{}{${}^*$ Supported in part by the NSF-grant DMS-0070816 and the MSRI}
\footnote{}{${}^\dagger$ Supported in part by the NSF-grant DMS-0097314 and the MSRI}

\def\date{October 29, 2001}
\def\Box #1 { \msk\par\nin 
\centerline{
\vbox{\offinterlineskip
\hrule
\hbox{\vrule\strut\hskip1ex\hfil{\smc#1}\hfill\hskip1ex}
\hrule}\vrule}\msk }

\def\address
{Simon Gindikin

Department of Mathematics

Rutgers University

New Brunswick, NJ 08903

USA

{\tt gindikin@math.rutgers.edu}

}

\def\addresstwo
{Bernhard Kr\"otz

The Ohio State University 

Department of Mathematics 

231 West 18th Avenue 

Columbus, OH 43210--1174 

USA

{\tt kroetz@math.ohio-state.edu}
}

\firstpage 

\subheadline{Abstract}

In this paper we define a distinguished boundary for  the 
complex crowns $\Xi\subeq G_\C /K_\C$ 
of non-compact Riemannian symmetric spaces $G/K$. The basic result is that 
affine symmetric spaces of $G$ can appear as a component of this boundary 
if and only if they are non-compactly causal symmetric spaces.

\sectionheadline{Introduction}

Let $X=G/K$ be a semisimple non-compact Riemannian symmetric space.  
We may assume that $G$ is semisimple with finite center. We write $G=NAK$ 
for an Iwasawa decomposition of $G$. By our assumption,  $G$ 
sits in its universal complexification $G_\C$ and so 
$X\subeq X_\C\:=G_\C/K_\C$. Note that 
$G$ does not act properly on $X_\C$.
The {\it complex crown} $\Xi\subeq X_\C$ of $X$ (cf.\ [Gi98]) 
was first considered  in [AkGi90]. It can be defined by 
$$\Xi\:=G\exp(i\Omega)K_\C/ K_\C$$
where $\Omega$ is a polyhedral convex domain in $\a\:={\rm Lie}(A)$
defined by 
$$\Omega\:=\{X\in\a\: (\forall\alpha\in\Sigma)\ |\alpha(X)|<{\pi\over 2}\}.$$
Here $\Sigma$ denotes the restricted root system with respect to $\a$.
Note that $G$ acts properly on $\Xi$ (cf.\ [AkGi90]). 

\par The domain $\Xi$ is interesting in several ways.  It accumulates
many crucial geometrical and analytical properties of $X$. For example, 
it was shown in [GiMa01] that  $\Xi$ 
parametrizes the compact cycles in complex flag domains.  
{}From the point of harmonic analysis the domain $\Xi$ is universal 
in the sense that all eigenfunctions on $X$ for the algebra 
$\D(X)$ of $G$-invariant differential operators extend 
holomorphically to $\Xi$ (cf.\ [KrSt01a]). It was conjectured in  [AkGi90] 
that the domain $\Xi$ is Stein; it was proved for different cases by 
different authors during the last year (cf.\ [GiKr01] for 
references and a more detailed account on the domain $\Xi$). 

\par In the case of groups $G$ of Hermitian type, $X$ is 
is a Hermitian symmetric space and $\Xi$ is $G$-equivariantly 
biholomorphic to $X\times \oline X$, where $\oline X$ refers to 
$X$ equipped with the opposite complex structure (cf.\ [BHH01], [GiMa01]
or [KrSt01b]). This example can be seen in a more 
general framework. Suppose that 
$X$ is the real form a Hermitian symmetric space. 
This means that there is a Hermitian 
group $S$ containing $G$ with maximal compact subgroup $U\supeq K$
such that the inclusion  $G/K\into S/U$ realizes $X$ as 
a totally real submanifold of $S/U$. For example, if $G$ is 
Hermitian we take $S=G\times G$ and $U=K\times K$. The existence 
of $S/U$ is guaranteed in many cases. In particular, $S/U$ exists 
for all classical groups $G$ in the sense that we have to replace 
$G$ sometimes with $G\times\R$ (e.g. $\Sl(n,\R)$ with $\Gl(n,\R)$). 
By the results of [BHH01] or [KrSt01b] there exists 
a generic $G$-invariant subdomain $\Xi_0$ of $\Xi$ which 
is $G$-biholomorphic to $S/U$. Moreover, 
equality $\Xi_0=\Xi$ holds if and only 
if $\Sigma$ is of type $C_n$ or $BC_n$ for $n\geq 2$ 
or $G=\SO(1,n)$ with $S=\SO(2,n)$. 
For the cases where $\Xi_0\subsetneq \Xi$, for example 
if $G=\SO(p,q)$ ($p,q> 2$) the geometric structure of $\Xi$
becomes very complicated and especially complicated
can be the boundary of $\Xi$.

\par In this paper we start to investigate the boundaries of 
$\Xi$. Write $\partial\Xi$ for the boundary of $\Xi$ in $X_\C$.
Apparently, the boundary will be a union of $G$-orbits, but the 
stratification of these orbits can be intricate. 
Firstly, not all orbits on the boundary intersect $\exp(i\a)K_\C/ K_\C$ --
but in this paper we focus only on such orbits.
Moreover, we are only interested in very special orbits of such type 
which are in some sense minimal. Let us define 
 the {\it distinguished boundary} $\partial_d\Xi$
of $\Xi$ as the union of the following $G$-orbits in $\partial\Xi$:
$$\partial_d\Xi\:=G\exp(i\partial_e\Omega)K_\C/ K_\C, $$
where $\partial_e\Omega$ the set of extreme points in the polyhedral
compact convex set $\oline{\Omega}$. Note that  
$$\partial_e\Omega={\cal W}(Y_1)\amalg\ldots\amalg{\cal W}(Y_n)$$
is a finite union of Weyl group orbits, where ${\cal W}$ is 
the Weyl group of $\Sigma$.

\par The distinguished boundary is a geometrically complicated object. 
Usually it is a disconnected set. Nevertheless, we show that it is minimal 
from some analytical points of view and that it features 
properties expected from a Shilov boundary. In particular, 

\msk\nin{\bf Theorem A.} {\it Write ${\cal A}(\Xi)$ for the algebra
of bounded holomorphic functions which continuously extend to the 
boundary of $\Xi$. Then 
$$(\forall f\in {\cal A}(\Xi))\qquad \sup_{z\in \Xi} |f(z)|
=\sup_{z\in \partial_d \Xi} |f(z)|.$$
Further $\partial_d\Xi\subeq \partial\Xi$ is minimal in a certain 
plurisubharmonic sense as explained in Section 1.}
\msk 

In above mentioned cases when  $\Xi_0=\Xi$ the situation is simpler. 
We can realize the Hermitian domain $D=S/U$ via the  Cartan embedding in 
the compact Hermitian symmetric $Y$ space dual to $S/U$. Here 
it has a compact boundary  and a compact Shilov boundary $\partial _c D$. 
Further, in this situation 
the Stein manifold $X_\C$ can be realized as a Zariski open part 
of $Y$ which will contain $\Xi$ biholomorphically equivalent to $D$ and 
$\partial_d \Xi$ will be only a Zariski open part of the compact manifold 
$\partial_c D$.  Let us emphasize that the Shilov boundary of $D$ 
essentially depends on the realization of $D$. 

\par We describe now the distinguished  boundary $\partial_d\Xi$ in more detail. 
Write $z_j'\:=\exp(iY_j)K_\C\in\partial_d\Xi$ for all $1\leq j\leq n$. Denote
by $H_j$ the isotropy 
subgroup of $G$ in $z_j'$. Then it follows from the definition of the 
distinguished boundary that we have a $G$-isomorphism
$$\partial_d\Xi\simeq G/H_1\amalg\ldots\amalg G/H_n.$$

\par In [Gi98] it was conjectured that non-compactly 
causal symmetric spaces appear in the ``Shilov boundary'' 
of the complex crowns $\Xi$. We establish this conjecture (in 
a more exact form);  namely we prove: 

\msk\nin{\bf Theorem B.} {\it If one of the components 
$G/H_j$ in $\partial_d\Xi$ is a symmetric space, then 
it is a non-compactly causal symmetric space. Moreover, 
every non-compactly causal symmetric space occurs as 
a component of the distinguished  boundary of some complex crown $\Xi$.}

\msk Let us say a few words about the motivation of this conjecture. 
On Riemannian symmetric spaces we have an elliptic analysis and on
non-compactly causal symmetric spaces we have a hyperbolic analysis. 
It is known in mathematical physics that in many important  cases elliptic 
and hyperbolic theories can be ``connected'' through complex domains
(Laplacians and wave equations, Euclidean and Minkowski field theories etc). 
Theorem B implies a connection of  Riemannian  and non-compactly 
causal symmetric spaces through the complex crowns $\Xi$. It shows 
that the phenomenon described above has a non-trivial generalization 
to symmetric spaces.

\par Write $\g\:={\rm Lie}(G)$ for the Lie algebra $G$. 
We call $\g$ {\it non-compactly causal} if $\g$ admits
an involution such that $(\g,\tau)$ is a non-compactly 
causal symmetric Lie algebra. Theorem B then tells us that 
$\g$ has to be non-compactly causal in order for 
$\partial_d\Xi$ to contain symmetric spaces. The class of non-compactly 
causal Lie algebras is very rich; for example it contains all 
classical Lie algebras except $\su(p,q)$ for $p\neq q$, $\sp(p,q)$
for $p\neq q$ and $\so^*(2n)$ for $n$ odd (cf.\ [Hi\'Ol96] or 
the table in Remark 3.3(a) below). 

\par In Theorem 3.24 we give a complete classification of 
the distinguished boundary $\partial_d\Xi$ for $\g$ non-compactly causal. 
It turns out that $\partial_d\Xi$ is a union of non-compactly 
causal symmetric spcaes except for
$\g=\so(p,q)$ with $3\leq p<q$ or $\g=\so(2n+3,\C)$,
where $\partial_d\Xi\cong G/H_1\amalg G/H_2$ with 
$G/H_1$ non-compactly causal and $H_2$ a non-symmetric subgroup of $G$.
Another consequence of the classification is the following 
result: 

\msk\nin {\bf Theorem C.} {\it Let $\g$ be a non-compactly causal Lie algebra
and $\partial_d\Xi\simeq\coprod_{j=1}^n G/H_j$ the decomposition 
of $\partial_d\Xi$ into $G$-orbits.  
Then a boundary component  $G/H_j$ of $\partial_d\Xi$
is totally real if and only if $G/H_j$ is symmetric.}

\msk  In Sections 4 and 5 we compare the dinstinguished boundary 
of $\Xi$ with the distinguished boundary of $\Xi_0$  for 
all those  groups $G$ for which $\Xi_0$ exists and 
$\Xi_0\subsetneq \Xi$. These are precisely the structure groups
of Euclidean Jordan algebras and the special orthogonal groups
$G=\SO_e(p,q)$ and $G=\SO(n,\C)$. 

\ssk This paper serves as the geometric foundation 
for the forthcoming work of the authors with Gestur \'Olafsson 
towards the definition of a Hardy space on $\Xi$ which 
realizes the most-continuous part $L^2(G/H)_{\rm mc}$ of $L^2(G/H)$
for a non-compactly causal symmetric space $G/H$ (cf.\ [GK\'O01]).

\msk It is our pleasure to thank the MSRI, Berkeley,  for 
its hospitality during the {\it Integral geometry program}
where this work was accomplished. We are grateful to Gestur \'Olafsson 
for going over the manuscript and his worthy suggestions.

\sectionheadline{1. Notation} 

Let $G$ be a semisimple Lie group sitting inside its universal  
complexification $G_\C$. 
We denote by $\g$ and $\g_\C$ the Lie algebras of $G$ and $G_\C$, respectively. 
Let $K<G$ be a maximal compact subgroup and $\k$ its Lie algebra. 

\par   Let $\g=\k\oplus\p$ be  the Cartan decomposition attached to $\k$. 
Take  $\a\subeq 
\p$ a maximal abelian subspace and let $\Sigma=\Sigma(\g, \a)\subeq \a^*$ be the 
corresponding root system. Related to this root system is the root space 
decomposition according to the simultaneous eigenvalues of $\ad (H), H\in \a:$

$$\g=\a\oplus \m \oplus\bigoplus_{\alpha\in \Sigma} \g^\alpha,$$ 
here $\m=\z_\k(\a)$ and $\g^\alpha =\{ X\in \g\: (\forall H\in \a) \ 
[H, X]=\alpha(H) X\}$. For the choice of a positive system $\Sigma^+\subeq 
\Sigma$ 
one obtains the nilpotent Lie algebra $\n =\bigoplus_{\alpha\in 
\Sigma^+}\g^\alpha$. Then one has the Iwasawa decomposition on the Lie algebra level 
$$\g=\n\oplus \a \oplus \k.$$ 
We write $A$, $N$ for the analytic subgroups of $G$ corresponding to 
$\a$ and $\n$.  
For these choices one has for $G$ the Iwasawa decomposition, namely, 
the multiplication map 

$$N\times A\times K\to G, \ \ (n, a, k)\mapsto nak$$
In particular, every element 
$g\in G$ can be written uniquely as $g=n(g) a(g) \kappa(g)$ with each of the 
maps $\kappa(g)\in K$, $a(g)\in A$, $n(g)\in N$ depending analytically on $g\in 
G$. 
The last piece of structure theory we shall recall is the little Weyl group. We 
denote by ${\cal W} =N_K(\a)/ Z_K(\a)$ the {\it Weyl group} of 
$\Sigma(\a,\g)$. 

Next we define a domain using the restricted roots. We set
$$\Omega=\{ X\in\a\:(\forall \alpha\in \Sigma)\  |\alpha(X)|<{\pi\over2}\}.$$
Clearly, $\Omega$ is convex and ${\cal W}$-invariant. 

\ssk If $L$ is a Lie group, then we write $L_0$ for the 
connected component containing $\1$. If $L$ is a group and 
$\sigma$ is an involution on $L$, 
then we write $L^\sigma$ for the $\sigma$-fixed points in $L$. 
Similarily, if $\l$ is a Lie algebra and $\sigma$ 
an involution on $\l$, then we write $\l^\sigma$ for the subspace 
of $\l$ which is pointwise fixed under $\sigma$.

\subheadline{The domain $\Xi$}
Write $A_\C$, $N_\C$ and $K_\C$ for the complexifications of the 
groups $A$, $N$ and $K$ realized in $G_\C$.  
Define the domain 
$$\Xi\:= G\exp(i\Omega) K_\C/ K_\C$$
which will be in the center of attention throughout this paper.  
Write $\partial\Xi$ for the topological boundary of $\Xi$
in $G_\C/ K_\C$. Recall from [AkGi90] that the piece of the boundary 
of $\Xi$ which intersects 
$\exp(i\a)K_\C/K_\C$ is given through 
$$\partial_a\Xi\:=G\exp(i\partial\Omega)K_\C/ K_\C.$$
Note the following properties of $\Xi$: 

\msk 
\item{${\bf\cdot}$} $\Xi$ is open in $G_\C/ K_\C$ (cf.\ [AkGi90]). 
\item{${\bf\cdot}$} $\Xi$ is connected and $G$-invariant. 
\item{${\bf\cdot}$} $G$ acts properly on $\Xi$ (cf.\ [AkGi90]).  
\item{${\bf\cdot}$} The domain $\Xi$ is Stein (special cases  of
[BHH01, Th.\ 10] or [GiKr01, Th.\ 3.4])
\item{${\bf\cdot}$} One has $\Xi\subeq N_\C A\exp(i\Omega)K_\C/ K_\C$
(special case of [GiKr01, Lemma 3.1])
\msk

\sectionheadline{2. The distinguished  boundary of $\Xi$}

Write $\oline\Xi$ for the closure of $\Xi$ in $G_\C/ K_\C$. 
Note that $\oline \Xi =\Xi\amalg \partial \Xi$. Since 
$\oline\Xi$ is not compact, there is no standard definition of a 
Shilov boundary for $\Xi$. However, it is possible to define a natural 
distinguished boundary $\partial_d\Xi$ which features many properties 
of a Shilov boundary as we  will explain below.  

\par We write ${\cal A}(\Xi)$ for the algebra of bounded 
continuous functions on $\oline \Xi$ 
which are holomorphic when restricted to $\Xi$:

$${\cal A}(\Xi)\:=\{ f\in C(\oline \Xi)\: f\res_{\Xi}\in {\cal O}(\Xi), 
\ \|f\|\:=\sup_{z\in \Xi} |f(z)|<\infty\}.$$

It is easy to check that ${\cal A}(\Xi)$ equipped with the supremum 
norm is a commutative Banach algebra with identity. 
\par There is a natural action of $G$ on ${\cal A}(\Xi)$
by left translation in the arguments: 

$$G\times {\cal A}(\Xi)\to {\cal A}(\Xi), \ \ (g,f)\mapsto \lambda_g(f);\ 
\lambda_g(f)(z)=f(g^{-1}z).\leqno(2.1)$$

\msk The domain $\Xi$ is special in the sense that it admits 
a maximal abelian flat subdomain:

$$T_{\Omega}=A\exp(i\Omega)\subeq A_\C.$$
Write $\oline {T_{\Omega}}= A\exp(i\oline{\Omega})$ for the closure of 
$T_{\Omega}$ in $A_\C$. 
Note that $T_{\Omega}$ is biholomorphic to the tube domain 
$\a+i\Omega\subeq \a_\C$ over $\Omega$ via the mapping 

$$\a+i\Omega\to T_{\Omega}, \ \ Z\mapsto \exp(Z).$$
Observe that $\oline{\Omega}$ is a compact polyhedron in $\a$. 
Denote by $\partial_e \Omega$ the extreme points of 
$\oline{\Omega}$. Notice that $\partial_e \Omega$ is a finite 
${\cal W}$-invariant set, hence a finite union of 
${\cal W}$-orbits:
$$\partial_e\Omega={\cal W}(Y_1)\amalg\ldots\amalg
 {\cal W}(Y_n).$$

\Remark 2.1. If the restricted root system $\Sigma$ is of type 
$C_n$, $BC_n$ then $\partial_e \Omega$ consists 
of a single ${\cal W}$-orbit (see our  discussion in Section 3). \qed 

Set 
$$\partial_d T_{\Omega}=A\exp(i\partial_e \Omega).$$

\par Similarily as before we define the Banach algebra ${\cal A}(T_
{\Omega})$. 
Then we have the following elementary lemma:

\Lemma 2.2. For all $f\in {\cal A}(T_{\Omega})$ we have 

$$\sup_{z\in T_{\Omega}} |f(z)|=\sup_{z\in \partial_d T_{\Omega}} |f(z)|.$$

\Proof. Recall that $T_{\Omega}$ is isomorphic to the tube domain 
$\a+i\Omega$ with polyhedral bounded base $\Omega$. For 
$\dim\a=1$ this lemma is just the Phragmen-Lindel\"of 
theorem. An easy iteration of this argument gives a proof 
in the general case. 
\qed

Note that the mapping 

$$\oline{T_{\Omega}}\into \oline \Xi, \ \ a\mapsto aK_\C$$
defines a continuous embedding which is holomorphic 
when restricted to $T_{\Omega}$. If $x_0\:=K_\C$ denotes 
the base point in $G_\C/K_\C$, then 

$$GT_{\Omega}(x_0)=\Xi, \quad G\oline{T_{\Omega}}(x_0)\subeq 
\oline \Xi, \quad \hbox 
{and}\quad G\partial T_{\Omega}(x_0)\subeq \partial_d\Xi.\leqno(2.2)$$

Finally we define 

$$\partial_d\Xi= G\exp(i\partial_e\Omega)K_\C/K_\C\subeq \partial \Xi.$$

We call  $\partial_d\Xi$ the {\it distinguished boundary of $\Xi$
in $G_\C/K_\C$}. 
This notion is justified by the following result:

\Theorem 2.3. For all  $f\in {\cal A}(\Xi)$ we have 
$$\sup_{z\in \Xi} |f(z)|=\sup_{z\in \partial_d \Xi} |f(z)|.$$

\Proof.  Recall from (2.1) the natural action of $G$ on ${\cal A}(\Xi)$. 
With Lemma 2.2 and (2.2) we therefore obtain that

$$\eqalign{\sup_{z\in \Xi} |f(z)|&
=\sup_{g\in G}\sup_{z\in T_{\Omega}} |f(g^{-1}z(x_0))|
=\sup_{g\in G}\sup_{z\in T_{\Omega}} |\lambda_g(f)(z(x_0))|\cr 
& =\sup_{g\in G}\sup_{z\in \partial_d T_{\Omega}} |\lambda_g(f)(z(x_0))|
=\sup_{g\in G}\sup_{z\in \partial_d T_{\Omega}} |f(g^{-1}z(x_0))|\cr 
&=\sup_{z\in \partial_d\Xi} |f(z)|.\cr}$$
This proves the theorem.\qed 

\Remark 2.4. If $\partial_d\Xi$ is connected, i.e., if 
$\partial_e\Omega$ consists of a single 
${\cal W}$-orbit (cf.\ Remark 2.1), then $\partial_d\Xi$ is 
in fact minimal in the sense of Theorem 2.3. Hence it makes
perfect sense to call $\partial_d\Xi$ the {\it Shilov boundary}
of $\Xi$ in this case.\qed

\subheadline{Minimality of the distinguished boundary}

The whole boundary $\partial\Xi$ is a rather complicated 
stratified $G$-space. Also  the piece 
$\partial_a\Xi$ is only a very small subset in $\partial\Xi$. 

\par Because of the complicated structure of $\partial\Xi$ we 
were not able to prove that the distinguished boundary $\partial_d\Xi$ is minimal 
in $\partial\Xi$ in the sense of Theorem 2.3. However, we can obtain 
the minimality of $\partial_d\Xi$ if we replace 
the algebra ${\cal A}(\Xi)$ by the Banach algebra 
${\cal P}(\Xi)$ of bounded continuous plurisubharmonic  functions on $\Xi$:

$${\cal P}(\Xi)\:=\{ f\in C(\Xi)\: f \ \hbox{plurisubharmonic}, \ \  
 \|f\|\:=\sup_{z\in \Xi} |f(z)|<\infty\}.$$

For $0\leq t<1$ define  
$$\partial_{d,t}\Xi= G\exp(it\partial_e\Omega)K_\C/K_\C.$$
and note that $\lim_{t\to 1}\partial_{d,t} \Xi =\partial_d\Xi$ 
set-theoretically.  

\Theorem 2.4. The following assertions hold: 

\item{(i)} For all $f\in {\cal P}(\Xi)$ we have 
$$\sup_{z\in \Xi} |f(z)|=\lim_{t\to 1}\big(\sup_{z\in \partial_{d,t} \Xi} 
|f(z)|\big)\ .$$
\item {(ii)} $\partial_d\Xi$ is minimal in 
$\partial_a\Xi$ with respect to the property in (i). 

\Proof. (i) This is proved in the same way as Theorem 2.3. 
\par\nin (ii) Recall the partition 
$\partial_e\Omega={\cal W}(Y_1)\amalg\ldots\amalg{\cal W}(Y_n)$. 
Further set $\partial_{d,t,j}\Xi\:=G\exp(it  Y_j)K_\C/ K_\C$
for every $1\leq j\leq n$ and $0\leq t<1$.  
\par Suppose that $\partial_d\Xi$ is not minimal with 
respect to the property in (i). Since a minimal set is necessarily 
$G$-invariant, we hence find an $Y_j$ such that 
$$\lim_{t\to 1}\big(\sup_{z\in \partial_{d,t,j}\Xi} |f(z)|\big) \leq 
\lim_{t\to 1}\max_{k\neq j} \big(\sup_{z\in \partial_{d,t,k}\Xi} |f(z)|\big)
\leqno(2.3)$$
for all $f\in {\cal P}(\Xi)$.

\par To conclude the proof of (ii) we have to recall some facts from 
[GiKr01, Sect.\ 3]. 
There exists a unique 
holomorphic surjective mapping 
$$\Xi\to T_{\Omega}\subeq A_\C,  \ \ z\mapsto a(z)$$
such that $z\in N_\C a(z)K_\C/ K_\C$. Furthermore if 
$z=g\exp(iX)K_\C$ for $g\in G$ and $X\in\Omega$, then 
$$\Im \log a(z)\subeq \conv({\cal W}(X)).$$
Recall that $T_\Omega$ is biholomorphic 
to $\a+i\Omega$. Thus we can find an $A$-invariant
$F\in {\cal P}(T_\Omega)$ which peaks 
at $\exp(iY_j)$. Now the function 
$$f\: \Xi\to \C,  \ \ z\mapsto F(a(z))$$
is a well defined element in ${\cal P}(\Xi)$. By construction 
$f$ does not satisfy (2.3), concluding the proof of (ii).\qed

\sectionheadline{3. Determination of  the distinguished boundary}

In this section we will show that 
the distinguished boundary $\partial_d\Xi$ can only be a union of 
symmetric spaces under the geometric assumption that 
$G$ is a non-compactly causal group. Further we will 
determine the distinguished boundary for all non-compactly 
causal groups.

\par We first have to recall some facts on causal symmetric 
Lie algebras.

\subheadline{Causal symmetric Lie algebras}

\par The standard reference for the facts collected below is the 
book [Hi\'Ol96].

\par As before $\g$ denotes a semisimple real Lie algebra and 
$\theta$ a Cartan involution on $\g$ with Cartan decomposition $\g=\k\oplus 
\p$. Let $\tau\: \s\to\s$ be an involution on $\g$ of which we may 
assume to commute with $\theta$. Write $\g=\h\oplus\q$ for the 
$\tau$-eigenspace decomposition 
corresponding to the $\tau$-eigenvalues $+1$ and $-1$. 

\par The symmetric Lie algebra $(\g,\tau)$ is called {\it irreducible} if the 
only $\tau$-invariant ideals of $\g$ are $\{0\}$ and $\g$. 

\Remark 3.1. If $(\g,\tau)$ is irreducible, then $\g$ is 
either simple or $\g=\h\oplus\h$ with $\h$ simple and 
$\tau$ is the flip $\tau(X,Y)=(Y,X)$ (``group case'').\qed 

\par If $(\g,\tau)$ is a symmetric Lie algebra, then we define its 
{\it c-dual} by $\g^c=\h\oplus i\q$. Call the restriction of 
the complex linear 
extension of $\tau$ to $\g^c$ also by $\tau$.
Then $(\g^c,\tau)$ is called the symmetric Lie algebra {\it c-dual}
to $(\g,\tau)$.

\Definition 3.2. Let $(\g,\tau)$ be an 
irreducible semisimple symmetric Lie algebra. 
\par\nin (a) We call $(\g,\tau)$ 
{\it compactly causal}
if $\z(\k)\cap\q\neq \{0\}$.
\par\nin (b) We call $(\g,\tau)$ 
{\it non-compactly causal} if its $c$-dual is compactly causal, 
or, equivalently, if there exists an $0\neq Y_0\in \q\cap\p$ which 
is fixed under $\h\cap\k$.  
\par\nin (c) We call $(\g,\tau)$ of 
{\it Cayley type} if it is both compactly and non-compactly 
causal. 
\qed 

\Remark 3.3. (a) Non-compactly causal 
symmetric Lie algebras are classified. The 
complete list is as follows (cf.\ [Hi\'Ol96, Th.\ 3.2.8]): 

$$\vbox{\tabskip=0pt\offinterlineskip
\def\tablerule{\noalign{\hrule}}
\halign{\strut#&\vrule#\tabskip=1em plus2em&
\hfil#\hfil&\vrule#&\hfil#\hfil &\vrule#\tabskip=0pt\cr\tablerule
&&\omit\hidewidth $\g$\hidewidth&& 
\omit\hidewidth $\h$ \hidewidth & \cr\tablerule
&& $\sp(n,\R)$ && $\gl(n,\R)$ & \cr\tablerule
&& $\su(n,n)$ && $\sL(n,\C)\times\R$ &  \cr\tablerule
&& $\so^*(4n)$  && $\sL(n,\H)\times\R$ & \cr\tablerule
&& $\so(p,q)$  && $\so(1,p-1)\times \so(1,q-1)$ & \cr\tablerule 
&& $\so(n,n)$  && $\so(n,\C)$ & \cr\tablerule 
&& $\sp(n,n)$  && $\sp(n,\C)$ & \cr\tablerule 
&& $\sL(n,\R)$ && $\so(q,n-q)$ \ $(1\leq q< n)$ & \cr\tablerule
&& $\sL(n,\H)$ && $\sp(q,n-q)$ \ $(1\leq q< n)$ & \cr\tablerule
&& $\so(2n,\C)$ && $\so^*(2n)$ & \cr\tablerule
&& $\so(n+2,\C)$ && $\so(2,n)$ & \cr\tablerule 
&& $\sp(n,\C)$ && $\sp(n,\R)$ & \cr\tablerule 
&& $\sL(n,\C)$ && $\su(q,n-q)$\ $(1\leq q< n)$ & \cr\tablerule
&& $\e_{6(6)}$ && $\sp(2,2)$ & \cr\tablerule
&& $\e_{6(-26)}$ && $\f_{4(-20)} $ & \cr\tablerule
&& $\e_6$ && $\e_{6(-14)} $ & \cr\tablerule
&& $\e_{7(7)}$ && $\su^*(8)$ & \cr\tablerule
&& $\e_{7(-25)}$ && $\e_{6(-26)}\times\R$ & \cr\tablerule
&& $\e_7$ && $\e_{7(-25)} $ & \cr\tablerule
}}$$

\par\nin (b) If $(\g,\tau)$ is of Cayley type, then $\g$ has to be simple 
hermitian and of tube type. Conversely, if $\g$ is simple 
hermitian and of tube type, then there exists up to conjugation
only one involution $\tau$ on $\g$ (the square of the Cayley transform)
turning $(\g,\tau)$ into 
a Cayley type symmetric Lie algebra (cf.\ [Hi\'Ol96]).
The Cayley type spaces are the following:
$$\vbox{\tabskip=0pt\offinterlineskip
\def\tablerule{\noalign{\hrule}}
\halign{\strut#&\vrule#\tabskip=1em plus2em&
\hfil#\hfil&\vrule#&\hfil#\hfil &\vrule#\tabskip=0pt\cr\tablerule
&&\omit\hidewidth $\g$\hidewidth&& 
\omit\hidewidth $\h$ \hidewidth & \cr\tablerule
&& $\sp(n,\R)$ && $\gl(n,\R)$ & \cr\tablerule
&& $\su(n,n)$ && $\sL(n,\C)\times\R$ &  \cr\tablerule
&& $\so^*(4n)$  && $\sL(n,\H)\times\R$ & \cr\tablerule
&& $\so(2,n)$ && $\so(1,n-1)\times\R$ & \cr\tablerule 
&& $\e_{7(-25)}$ && $\e_{6(-26)}\times\R$ & \cr\tablerule
}}$$
\qed 

We call a Lie algebra $\g$ {\it non-compactly causal} if there 
exists an involution $\tau$ on $\g$ turning $(\g,\tau)$ into 
a non-compactly causal symmetric Lie algebra. A Lie group 
$G$ is called {\it non-compactly causal} if its 
Lie algebra $\g$ is non-compactly causal.

\subheadline{General results on $G$-orbits in the distinguished boundary}

Let $\partial_d\Xi=G\exp(i\partial_e\Omega)K_\C/ K_\C$ be the 
distinguished boundary of $\Xi$. 
Let 
$$\partial_e\Omega={\cal W}(Y_1)\amalg\ldots\amalg
 {\cal W}(Y_n)$$
be the partition of Weyl group orbits with 
$Y_j\in \partial_e\Omega$. Set $z_j\:=\exp(iY_j)$
and $z_j'\:=zK_\C\in\partial_d\Xi$.  
Recall that 
$$\partial_d\Xi\simeq\coprod_{j=1}^n G/G_{z_j'}$$
where $G_{z_j'}$ denotes the isotropy subgroup of $G$ 
in $z_j'$. 

\par We also write $\theta$ for the holomorphic 
extension of the Cartan involution from $G$ to $G_\C$. 

\par If we assume that $G_\C^\theta=K_\C$ 
(this is always satisfied if $G_\C$ is simply connected), 
then  we can characterize $G_{z_j'}$ through the following 
equivalences: 

$$\eqalign{g\in G_{z_j'} &\iff z_j^{-1}g z_j\in K_\C\cr 
&\iff (z_j^{-1} g z_j) \theta(z_j^{-1} g z_j)^{-1}=\1\cr 
&\iff g z_j^2\theta(g)^{-1}=z_j^2\cr
&\iff z_j^2\theta(g)z_j^{-2}=g.\cr}\leqno(3.1)$$
Note that we always have $(G_\C^\theta)_0=K_\C$. 
If $G_\C^\theta\neq K_\C$, then the computation (3.1)
still gives us the connected component $(G_{z_j'})_0$
of $G_{z_j'}$, namely 
$$(G_{z_j'})_0=\{ g\in G\: z_j^2\theta(g)z_j^{-2}=g\}_0.$$

\par Summarizing our discussions from above we have 
proved:

\Lemma 3.4. Let $Y\in \partial_e\Omega$ and set $z\:=\exp(iY)$, 
$z'\:=zK_\C\in\partial_d\Xi$. Then the following 
assertions hold: 
\item{(i)} The connected component of 
the isotropy subgroup in $z'$ is given by 
$$(G_{z'})_0=\{ g\in G\: z^2\theta(g)z^{-2}=g\}_0.$$
\item{(ii)} If $G_\C$ is simply connected, then 
$$G_{z'}=\{ g\in G\: z^2\theta(g)z^{-2}=g\}.$$
\item{(iii)} If $G_\C$ is simply connected 
and if $\tau(g)=z^2\theta(g)z^{-2}$ is an involution on $G$, then 
$$g\in G_{z'}\iff g\in G^\tau.\qeddis

In general we cannot expect that the automorphism $\tau$ in Lemma 3.4(iii)
preserves  $G$. This can only happen under 
the geometric assumption that $G$ is a non-compactly 
causal group.

\Theorem 3.5. Let $Y\in \partial_e\Omega$ and $z\:=\exp(iY)$. 
Consider the automorphism 
$$\tau\: G_\C\to G_\C, \ \ g\mapsto z^2\theta(g)z^{-2}.$$
Then the following assertions hold: 
\item{(i)} If $\tau\:=d\tau(\1)$ preserves $\g$,
then $\tau$ is an involution and $(\g,\tau)$
is a non-compactly causal symmetric Lie algebra.  
\item{(ii)} If $\tau$ preserves $G$, then $G/G^\tau$ is 
non-compactly causal and we have a $G$-isomorphism
$$G/G^\tau\simeq G(z')\subeq \partial_d\Xi.$$

\Proof. (following a suggestion of Gestur \'Olafsson) (i) Suppose 
that $\tau$ defines an automorphism of $\g$. Note that 
$\tau(X)= e^{i2\ad Y}(X)$ for $X\in\g$. 
\par Since $Y\in\oline \Omega$ we have 
$\Spec \ad (Y)\subeq [-{\pi\over 2}, {\pi\over 2}]$. 
Further the fact that $\tau(X)\in\g$ for all $X\in\g$ means precisely 
that $\Spec(\ad Y)\subeq \Z {\pi\over 2}$. 
Hence we obtain that $\Spec(\ad Y)=\{-{\pi\over 2}, 0, {\pi\over 2}\}$
by the symmetry of the spectrum. 
In particular, $\tau$ is an involution on $\g$
and it follows from [Hi\'Ol96, Th.\ 3.2.4] that 
$(\g,\tau)$ is non-compactly causal.
\par\nin (ii) This follows from (i).\qed

In view of Theorem 3.5, we can only expect that symmetric 
spaces are contained in the distinguished 
boundary $\partial_d\Xi$ if $G$ is a non-compactly causal group. 
Before we start to classify the distinguished boundaries related to 
all such groups we first discuss the case where 
$G$ is an arbitrary complex semisimple Lie group.

\subheadline{The case where $G$ is complex}

We will recall some material on  complexifications of semisimple 
complex Lie algebras and certain complex homogeneous spaces.

\par Let $\g$ be a complex semisimple Lie algebra and 
$\k<\g$ a maximal compact subalgebra.   Then $\k$ is 
a compact real form of $\g$ and $\g=\k + i\k$ is 
a Cartan decomposition of $\g$. We write $X\mapsto\oline X$
for the complex conjugation of $\g$ with respect to the 
real form $\k$. 

\par Write $\g_\R$ for $\g$ 
considered as a real Lie algebra and define $\g_\C$ as the 
complexification of $\g_\R$. Denote by $J$ the multiplication 
with $i$ in $\g_\C$. Finally we write $\g^{\rm opp}$ for $\g$ 
equipped with the opposite complex structure. 

The mapping 
$$\g_\C \to \g\oplus\g^{\rm opp},\ \ X+JY\mapsto (X+iY, X-iY)$$
is an isomorphism of the complex Lie algebra $(\g_\C, J)$ with 
$\g\oplus\g^{\rm opp}$. 
Now 
$$\g^{\rm opp}\to\g, \ \ X\mapsto \oline X$$
establishes a complex Lie algebra isomorphism. Hence we get 
the comlex Lie algebra isomorphism 
$$\g_\C \to \g\oplus\g, \ \ X+JY\mapsto (X+iY, \oline X +i\oline Y).$$
Note the following realizations of subalgebras in $\g_\C$ inside 
of $\g\oplus\g$: 

$$\g\simeq\{ (X, \oline X)\: X\in\g\}$$
$$\k_\C \simeq\Delta(\g)\:=\{ (X,  X)\: X\in\g\}$$
$$\a_\C\simeq\{ (X, - X)\: X\in\a+i\a\subeq \g\}$$

Let now $G$ be a Lie group with Lie algebra $\g$. 
Then our observations from above imply the 
isomorphism 
$$G_\C/K_\C\simeq (G\times G)/ \Delta(G)$$
with $\Delta(G)$ the diagonal subgroup. 
Further we have the canonical isomorphism 

$$G\times G/ \Delta (G)\to G, \ \ (g,h)\Delta(G)\mapsto gh^{-1}.$$
Hence we have $G_\C/K_\C\simeq G$.
 Note that the morphism  of $A_\C$
into $G_\C/ K_\C=G$ is given by 

$$A_\C\to G, \ \ a\mapsto a^2,$$ 
the square mapping. Further the action of $G$ on $G_\C/K_\C=G$ is 
given by 

$$G\times G_\C/K_\C\to G_\C/K_\C, \ \ (g,x)\mapsto gx\oline g^{-1},$$  
where we identified $G_\C/K_\C$ with $G$.

The summary of our discussion is now: 

\Proposition 3.6. Assume that $G$ is a semisimple complex Lie group. 
Then $G_\C/ K_\C$ can be canonically identified with $G$. 
Moreover, if $Y\in \partial_e\Omega$, $z=\exp(iY)$ and $z'=zK_\C$, 
then the isotropy subgroup of $G$ in $z'$ is given by 
$$G_{z'}=\{g\in G\: z^2 \oline g z^{-2} =g \},$$
where $g\mapsto \oline g$ denotes the conjugation in $G$ with 
respect to the compact real form $K$ of $G$.\qed

From now on we will assume that $\g$ is a non-compactly causal
Lie algebra. Depending on the type of the restricted root system
$\Sigma$ we are 
are now going to determine the distinguished boundary of $\Xi$. 

\subheadline{The cases with restricted root system of type $C_n$}

Assume that $\g$ is a non-compactly causal Lie algebra
with restricted root system of type $C_n$. 
According to the list in Remark 3.3(a) this means that $\g$ is
one of the following list: 

$$ \sp(n,\R)\quad \sp(n,n)\quad \sp(n,\C)\quad 
\su(n,n)\quad \so^*(4n)\quad \so(2,n)\quad \e_{7(-25)}.$$
Note that for all these cases there exists up to conjugation only  
one involution $\tau$  on $\g$ which turns $(\g,\tau)$ into 
a non-compactly causal symmetric Lie algebra. 

\par On the group level we take $G$ to be the real form of a simply connected
Lie group $G_\C$ with Lie algebra $\g_\C$. 

\par The restricted root system $\Sigma$ is of type $C_n$, say 

$$\Sigma=\{ {1\over 2}(\pm\eps_i\pm\eps_j)\: 1\leq i,j\leq n\}\bs\{0\}.$$
Define $Y_j\in \a$ by $\eps_i(Y_j)={\pi\over 2}\delta_{ij}$. 
Then 
$$\Omega=\bigoplus_{j=1}^n ]-1, 1[ Y_j$$
and 
$$\partial_e\Omega=\{\pm Y_1\pm\ldots\pm Y_n\}.$$
Define 
$$Y_0\:=Y_1+\ldots+Y_n$$
and note that 
$$\partial_e\Omega={\cal W} (Y_0),$$
since ${\cal W}$ consists of all permutations and 
sign changes. 
Observe that 
$$\Spec(\ad Y_0)=\{ -{\pi\over 2}, 0, {\pi\over 2}\}.$$

\ssk We have  $\partial_d\Xi=G\exp(iY_0)K_\C/K_\C$.  
Set $z_0\:=\exp(iY_0)$. Then 
$\partial_d\Xi=G(z_0)$ and so $\partial_d\Xi\simeq G/G_{z_0'}$
with $z_0'=z_0K_\C\in \partial_d\Xi$.
Note that 
$$\tau\: G_\C\to G_\C, \ \ g\mapsto z_0^2 \theta(g)z_0^{-2}$$
restricts to an involution on $G$ since 
$\Spec(\ad Y_0)=\{ -{\pi\over 2}, 0, {\pi\over 2}\}$. 
Thus we obtain from Lemma 3.4 and Theorem 3.5 that:

\Proposition 3.7. Let $G$ be a non-compactly causal group sitting 
inside a simply connected complex Lie group $G_\C$. 
Suppose that the restricted root system $\Sigma$ is of type $C_n$. 
Then the distinguished boundary of 
$\Xi$ is $G$-isomorphic to the (up to isomorphism unique)
non-compactly causal symmetric space $G/H$ associated to $G$: 
$$\partial_d\Xi\simeq G/H.\qeddis

\subheadline{Examples with restricted  root system of type $A_n$}

\nin{\bf The example of $G=\Sl(n,\R)$.}
Here we will encounter the situation that the distinguished boundary 
is not connected and a union of different non-compactly 
causal symmetric spaces. The example features the situation well 
where the restricted root system is of type $A_n$.

\par We introduce some notation. Let $\g=\sL(n,\R)$ and 
choose $\k=\so(n,\R)$ as a maximal compact subalgebra of $\g$.
Write $\g=\k\oplus\p$ for the corresponding Cartan decomposition. 
For $1\leq j\leq n$ we introduce the diagonal 
matrix $e_j\:=(\delta_{ij})_{i,j}$. Then 

$$\a\:=\{\sum_{j=1}^nx_j e_j\: x_j\in\R, \sum_{j=1}^n x_j =0\}$$
is a maximal abelian subspace in $\p$. 
Define linear functionals $\eps_j$ on the diagonal matrices 
by $\eps_j(e_i)=\delta_{ij}$. Then 

$$\Sigma=\{\eps_i-\eps_j\: i\neq j\}$$
and ${\cal W}$ is the permutation group of $\{ e_1, \ldots, e_n\}$. 

\par For every $1\leq q\leq n-1$ we define elements $Y_q\in\a$
by 
$$Y_q\:={\pi\over 2} \big(\sum_{j=1}^q e_j - {q\over n}\sum_{j=1}^n e_j\big).
$$
Then it follows from [KrSt01a, Lemma 2.1] that 

$$\partial_e\Omega=\coprod_{q=1}^{n-1} {\cal W} (Y_q).\leqno(3.2)$$
The union in (3.2) disjoint. 
For every $1\leq q\leq n-1$ set $z_q\:=\exp(iY_q)$ and 
denote by $H_q<G$ the stabilizer in $G$ of the element $z_qK_\C\in 
\partial_d\Xi$. 
Then we get 

$$\partial_d\Xi\simeq \coprod_{q=1}^{n-1} G/H_q.$$

It remains to calculate the subgroups $H_q$. 
Note that 
$$z_q=e^{-i{\pi\over2}{q\over n}}\diag(i,\ldots, i, 1,\ldots, 1)$$  
with $q$-times $i$ on the diagonal. From (3.1) we obtain the 
following equivalences:

$$\eqalign{g\in H_q &\iff z_q^{-1}g z_q\in K_\C\cr 
&\iff (z_q^{-1} g z_q)(z_q^{-1} g z_q)^t=\1\cr 
&\iff g z_q^2g^t=z_q^2\cr 
&\iff g \diag(-1,\ldots,-1,1,\ldots,1)g^t =\diag(-1,\ldots,-1,1,\ldots 1)\cr 
&\iff g\in \SO(q,n-q).\cr}$$
Alltogether we thus have shown:

\Proposition 3.8. For $G=\Sl(n,\R)$ we have the following 
$G$-isomorphism of the distinguished boundary:   

$$\partial_d\Xi\simeq \coprod_{q=1}^{n-1} \Sl(n,\R)/\SO(q,n-q).\qeddis

\msk 
\nin{\bf The example of $G=\Sl(n,\C)$.}
This is very similar to the case of $G=\Sl(n,\R)$. We will use 
our results on complex groups. Let $G=\Sl(n,\C)$ and 
$K=\SU(n)$. As before we can take $\a$ to be the diagonal 
matrices in $\g$. Since the root system is the same as for 
$\g=\sL(n,\R)$, we conclude that 

$$\partial_d\Xi\simeq \coprod_{q=1}^{n-1} G_{z_q'}$$
with $z_q'=\exp(iY_q)K_\C$ and $Y_q$ as before. 
\par Now we are going to use the identification of $G_\C/K_\C$
with $G$ as explained before Proposition 3.6. Then 
$z_q'\in G_\C/K_\C$ becomes identified with the 
point $z_q^2\:=\exp(i2Y_q)$ in $G$. We have 

$$z_q^2=c\diag(1, \ldots, 1, -1, \ldots, -1)$$
with $c$ a constant and $q$ times $+1$ on the diagonal. 
Further $G$ acts on $G_\C/K_\C=G$ by $g(x)=gx\oline g^t$
where  $\oline g$ denotes the complex conjugation with respect 
to $\Sl(n,\R)$. Hence Proposition 3.6 implies 
that $G_{z_q'}\simeq \SU(q,n-p)$ and we have proved:

\Proposition 3.9. For $G=\Sl(n,\C)$  we have the following 
$G$-isomorphism of the distinguished boundary:   

$$\partial_d\Xi\simeq \coprod_{q=1}^{n-1} \Sl(n,\C)/\SU(q,n-q).\qeddis

\msk
\nin{\bf The example of $G=\Sl(n,\H)$.}
Considerations very similar to the two other $A_n$-cases treated before 
yield the following result:

\Proposition  3.10. For $G=\Sl(n,\H)$ we have the following 
$G$-isomorphism of the distinguished boundary:

$$\partial_d\Xi\simeq \coprod_{q=1}^{n-1} \Sl(n,\H)/\Sp(q,n-q).\qeddis

\subheadline{Examples with restricted root systems of type $B_n$ or $D_n$}

The classical Lie algebras which have restricted root system of type 
$B_n$ or $D_n$ are $\so(p,q)$ and $\so(n,\C)$. Before 
we can treat these cases we first have to recall some facts 
of root systems of type $B_n$ and $D_n$. 

\msk\nin {\bf Root systems of  type $B_n$ and $D_n$.}
In the standard notation the root systems of type $B_n$  and 
$D_n$ are given by 

$$B_n=\{ \pm\eps_i\pm\eps_j\: 1\leq i\neq j\leq n\} \amalg\{ 
\pm\eps_i\: 1\leq i\leq n\}$$
and 
$$D_n=\{ \pm\eps_i\pm\eps_j\: 1\leq i\neq j\leq n\}.$$

As usual we write $e_1,\ldots, e_n$ for the dual basis 
to $\eps_1,\ldots, \eps_n$. 
If $\Sigma$ is a root systems, then we also write $\Omega=\Omega(\Sigma)$
to make the relation to $\Sigma$ clear whenever ambigious. 

\Lemma 3.11. For the root systems $B_n$ and $D_n$ the following 
assertions hold: 
\item{(i)} We have $\Omega(B_n)=\Omega(D_n)$. In particular, 
$\partial_e\Omega(B_n)=\partial_e\Omega(D_n)$.  
\item{(ii)} For $\Sigma=D_n$ and $n\geq 3$ we have 
$$\partial_e\Omega={\cal W}(Y_1)\amalg
{\cal W}(Y_2)\amalg
{\cal W}(Y_3)$$
a disjoint union with $Y_1={\pi\over 2}e_1$, $Y_2={\pi\over 4}(e_1+\ldots+e_n)$
and $Y_3={\pi\over 4}(e_1+\ldots+e_{n-1}-e_n)$.

\item{(iii)}For $\Sigma=B_2$ we have $\partial_e\Omega
={\cal W}(Y_1)$ and for $n\geq 3$ we have 
$$\partial_e\Omega={\cal W}(Y_1)\amalg
{\cal W}(Y_2)$$
a disjoint union with $Y_1$ and $Y_2$ as in (i).

\Proof. [KrSt01a, Sect.\ 2].\qed

\nin {\bf The example of $G=\SO(2n,\C)$ for $n\geq 4$.}
In this case the restricted root system is of type $D_n$. 
Hence Lemma 3.11(ii) implies that 
$$\partial_e\Omega={\cal W}(Y_1)\amalg{\cal W}(Y_2)\amalg
{\cal W}(Y_3)$$
consists of three Weyl group orbits. Note that $Y_2$ and $Y_3$ are 
conjugate under an outer isomorphism $\kappa$
induced from an outer isomorphism 
of the Dynkin diagram. 

\par As a maximal abelian subspace $\a\subeq\p$ we choose 
$$\a\:=\{ \pmatrix{\pmatrix{0 & it_1\cr -it_1 &0\cr} &  & \cr 
& \ddots &\cr 
& & \pmatrix{0 & it_n\cr -it_n &0\cr} \cr}\: t_1,\ldots,t_n\in\R\}.$$
If we set $z_j\:=\exp(iY_j)$ and 
$z_j'=z_j K_\C$ for $j=1,2,3$, then 
we obtain that 
$$z_1^2=\diag(-1,-1, 1,\ldots,1), \qquad 
z_2^2= \pmatrix{\pmatrix{0 & 1\cr -1 &0\cr} &  & \cr 
& \ddots &\cr 
& & \pmatrix{0 & 1\cr -1 &0\cr} \cr}$$
and $z_3^2=\kappa(z_2^2)$. 
From Proposition 3.6 we have 
$$G_{z_j'}=\{ g\in G\: z_j^2 \oline g z_j^{-2}=g\}$$ 
and so 
$$G_{z_1'}\simeq \SO(2, 2n-2)\quad\hbox{and}\quad 
G_{z_2'}\simeq G_{z_3'}\simeq \SO^*(2n).$$

Thus we have proved:

\Proposition 3.12. For $G=\SO(2n,\C)$ with $n\geq 3$ 
we have the following 
$G$-isomorphism of the distinguished boundary:
$$\partial_d\Xi\simeq \SO(2n,\C)/ \SO(2,2n-2)\amalg 
\SO(2n,\C)/ \SO^*(2n)\amalg \SO(2n,\C)/ \SO^*(2n).\qeddis

\msk\nin {\bf The example of $G=\SO(2n+1,\C)$.} 
In this case the restricted root system is of type $B_n$
and Lemma 3.11(iii) implies that 
$$\partial_e\Omega={\cal W}(Y_1)\amalg{\cal W}(Y_2)$$
where the second term in the union is not existent for $n=2$. 
As before we now obtain that:

\Proposition 3.13. For $G=\SO(2n+1,\C)$ we have the following 
$G$-isomorphism of the distinguished boundary:
\item{(i)} For $n\geq 3$ one has 
$$\partial_d\Xi\simeq \SO(2n+1,\C)/ \SO(2,2n-1)\amalg 
\SO(2n+1,\C)/ \SO^*(2n).$$
In particular, the boundary component 
$\SO(2n+1,\C)/ \SO^*(2n)$ is not a symmetric space. 
\item{(ii)} For $n=2$ one has 
$$\partial_d\Xi\simeq \SO(5,\C)/ \SO(2,3).\qeddis

\msk\nin {\bf Some general remarks on $\g=\so(p,q)$ for $0<p\leq q$.} 
An appropriate maximal abelian subspace 
of $\p$ in $\g=\so(p,q)$ is 
$$\a\:=\{ \pmatrix{ 0_{pp} & I_{t_1,\ldots, t_p}\cr 
I_{t_1,\ldots, t_p}^t & 0_{qq}\cr}\: t_1,\ldots,t_p\in \R\},$$ 
where 
$$I_{t_1,\ldots, t_p}=\pmatrix{ &  &  & t_1\cr 0_{p,p-q}& & \addots& \cr & t_p& & \cr}
\in M(p\times q;\R).$$

Let $n=p+q$. For all $1\leq i,j\leq n$ define $E_{ij}\in M_n(\R)$ by 
$E_{ij}=(\delta_{ki}\delta_{jl})_{k,l}$. Then 
$\a=\bigoplus_{j=1}^p \R e_j$ with 
$$e_j=E_{j, n+1-j}+ E_{n+1-j, j}.$$
Define $\eps_j\in \a^*$ by $\eps_j(e_k)\:=\delta_{jk}$. Then the
root system $\Sigma =\Sigma(\a,\g)$ is given by 
$$\Sigma=\cases{ \{ \pm\eps_i\pm\eps_j\: 1\leq i\neq j\leq p\} \amalg\{ 
\pm\eps_i\: 1\leq i\leq p\} & for $1<p< q$,\cr 
\{ \pm\eps_i\pm\eps_j\: 1\leq i\neq j\leq p\} & for $1<p=q$,\cr 
\{\pm \eps_1\} & for $p=1$.\cr}$$
The restricted  root system is hence of type $B_p$, $D_p$ or $A_1$.

\msk\nin {\bf The example of $G=\SO_e(n,n)$ for $n\geq 3$.}
Here the root system is of type $D_n$ and we have 
$$\partial_e\Omega={\cal W}(Y_1)\amalg{\cal W}(Y_2)\amalg
{\cal W}(Y_3).$$
In the coordinates introduced before we obtain that  

$$z_1^2=\diag(-1,1,\ldots, 1, -1), \qquad
z_2^2=\{ \pmatrix{ 0_{nn} & iI_{1,\ldots, 1}\cr 
-iI_{1,\ldots, 1} & 0_{nn}\cr}$$
and $z_3^2=\kappa(z_2^2)$. From Lemma 3.1 we thus obtain that 

$$(G_{z_1'})_0=\SO_e(1,n-1)\times\SO_e(1,n-1)\quad \hbox{and}\quad
(G_{z_2'})_0=\SO(n,\C).$$

Thus we have shown that:

\Proposition 3.14. For $G=\SO_e(n,n)$ and $n\geq 3$ we have the following 
$G$-isomorphism of the distinguished boundary: $\partial_d\Xi\simeq G/H_1\amalg 
G/H_2\amalg G/H_2$ with 
$$(H_1)_0=\SO_e(1,n-1)\times\SO_e(1,n-1)\quad (H_2)_0=\SO(n,\C).\qeddis

\msk\nin {\bf The example of $G=\SO(p,q)$ for $1\leq p<q$.}
Here the restricted  root system is of type $B_p$ for $p\geq 2$
and $A_1$ for $p=1$. Considerations very similar to 
as before lead us to the following result:

\Proposition 3.15. For $G=\SO_e(p,q)$ and $1\leq p<q$ the 
distinguished boundary $\partial_d\Xi$ is described as follows:

\item{(i)} For $p=1$ one has $\partial_d\Xi\simeq G/H$
with $H_0= \SO_e(1,q-1)$. 
\item{(ii)} For $p=2$ one has $\partial_d\Xi\simeq G/H$ with 
$H_0=\SO_e(1,q-1)\times\R.$
\item{(iii)} For $p\geq 3$ one has 
$\partial_d\Xi\simeq G/H_1\amalg G/H_2$ with 
$$(H_1)_0=\SO_e(1,p-1)\times\SO_e(1,q-1)\qquad \hbox{and}\qquad 
(H_2)_0=\SO(p,\C)\times \SO(q-p).$$
In particular, the boundary component $G/H_2$ is not a symmetric space. \qed

\subheadline{Exceptional cases}

We first have to recall some facts about the root systems 
$E_6$ and $E_7$. Our references herefore are the 
tables from [Kn96]. Write $e_1, \ldots, e_n$
for the standard basis in Euclidean space $\R^n$ and write $\la\cdot,\cdot
\ra$ for the standard inner product.

\msk\nin {\bf The system $E_6$.} Let 
$$V\:=\{ v\in \R^8\: \la v, e_6-e_7\ra=\la v, e_7+e_8\ra=0\}.$$
Then the root system $E_6$ can be defined by 

$$\Sigma\:=\{ \pm e_i\pm e_j\: 1\leq j<i\leq 5\}\amalg\{ {1\over 2}
\sum_{j=1}^8 (-1)^{n(j)} e_j\in V\: \sum_{j=1}^8 n(j) \ \hbox{even}\}.$$
A positive system is given by 
$$\Sigma^+=\{ e_i\pm e_j\: 1\leq j<i\leq 5\}\amalg\{ {1\over 2}
(e_8-e_7-e_6+ \sum_{j=1}^5 (-1)^{n(j)} e_j)
\: \sum_{j=1}^5 n(j) \ \hbox{even}\}.$$
Note that $\la \alpha,\alpha\ra =2$ for all $\alpha\in\Sigma$, i.e., 
all roots have the same length. 
A basis for $\Sigma^+$ is given by 

$$\Pi=\{\alpha_1, \ldots,\alpha_6\}=\{ {1\over 2}(e_8-e_7-e_6-e_5-e_4-e_3-
e_2+e_1), e_1 +e_2, e_2-e_1, e_3-e_2, e_4-e_3, e_5-e_4\}.$$
A simple calculation then shows that the fundamental weights 
are given by 

$$\eqalign{\omega_1 &={2\over 3}(e_8-e_7-e_6)\cr 
\omega_2 &= {1\over 2} (e_1+e_2+e_3+e_4+e_5)+ {1\over 2}(e_8-e_7-e_6)\cr 
\omega_3 &= {1\over 2} (-e_1 + e_2 +e_3 +e_4 +e_5) +{5\over 6}(e_8-e_7-e_6)\cr 
\omega_4 &=  (e_3 +e_4 +e_5) + (e_8-e_7-e_6)\cr 
\omega_5 &= (e_4+e_5) +{2\over 3} (e_8-e_7-e_6)\cr 
\omega_6 &=  e_5 + {1\over 3} (e_8-e_7-e_6).\cr}$$
Note that $\omega_1,\ldots,\omega_6$ is also the dual basis 
of $\Pi$ since $\la \alpha,\alpha\ra=2$ for all $\alpha\in \Sigma$. 

The highest root $\beta$ in $\Sigma^+$ can be writtem as 

$$\beta=\alpha_1+2\alpha_2+2\alpha_3+ 3\alpha_4+2\alpha_5+\alpha_6.$$

In the following we set $\a\:=V$ and we identify $\a$ with $\a^*$ via
the inner product on $\a$. 

If $E\subeq \a$ is a  closed convex set, then we write $\Ext(E)$ for the 
extreme points of $E$.

\Lemma 3.16. Let $\Sigma$ be an irreducible abstract root system in 
the Euclidean space $\a$. Let $\Pi=\{\alpha_1,\ldots,\alpha_n\}$
be a basis of $\Sigma$ and $\omega_1,\ldots,\omega_n$ its dual basis. Further 
let $\beta=\sum_{j=1}^n m_j \alpha_j$ be the highest root with 
respect to $\Pi$. Further write $C$ for the closed Weyl chamber with 
respect to $\Pi$. Then the following assertions hold: 

\item{(i)} $\Omega$ is ${\cal W}$-invariant. 
\item{(ii)} $\Ext(\oline {\Omega}) \cap C\subeq \Ext (\oline {\Omega}\cap C)$. 
\item{(iii)} We have the following inclusion
$$\Ext(\oline{\Omega})\subeq {\cal W} ({\pi\over2}\{ {\omega_1\over m_1},\ldots, 
{\omega_n \over m_n} \}).$$

\Proof. (i) and (ii) are trivial. 
\par\nin (iii) First note that (i) and (ii) imply that 
$$ \Ext(\oline{\Omega})\subeq {\cal W}(\Ext (\oline {\Omega})\cap C).$$
Thus it is sufficient to show that 
$$\Ext(\oline{\Omega})\cap C\subeq {\pi\over2}\{ {\omega_1\over m_1},\ldots, 
{\omega_n \over m_n}\}.$$
Note that 
$$C=\{ x\in \a\: (\forall j) \la x,\alpha_j\ra\geq 0\}=
\bigoplus_{j=1}^n \R^+ \omega_j.$$ 
Thus we get from $m_j\in \N$ for all $j$ that 
$$\partial \Omega\cap C=\{ x=\sum_{j=1}^n\lambda_j \omega_j\: \lambda_j\geq 0, 
\sum_{j=1}^n m_j\lambda_j={\pi\over2}\}.$$
From this the assertion follows. 
\qed 

For any $\alpha\in\a$ we denote by $s_\alpha$ the reflection 
with respect to the hyperplane $\alpha^\bot$.

\Lemma 3.17. For the root system $E_6$ we have 
$$\partial_e\Omega=\Ext(\oline \Omega)= {\cal W} ({\pi\over 2}\omega_1)\amalg {\cal W}(
{\pi\over 2} \omega_6)\ .$$

\Proof. First recall that $E_6$ admits a nontrivial 
outer automorphims $\kappa$, the reflection of the Dynkin diagram. 
Note that $\kappa$ leaves $\oline \Omega$ invariant, 
$\kappa(\alpha_1)=\alpha_6$
and $\kappa(\alpha_3)=\alpha_5$. 

\par In view of Lemma 3.16(iii), it is hence sufficient to 
show that ${\omega_2\over 2}$,  ${\omega_4\over 3}$, 
${\omega_5\over 2}$ are not extremal in $\oline {\Omega}$. 
Now observe that 
$$\eqalign{ {\omega_5\over 2}&={1\over 2}(e_4+e_5) +
{1\over 3}(e_8-e_7-e_6)= {1\over 2}(\omega_6)
+{1\over 2}(s_{e_4-e_5}\omega_6)\cr
{\omega_4\over 3}&={1\over 3}(e_3+e_4+e_5)+{1\over 3}(e_8-e_7-e_6)={2\over 3}({\omega_5\over 2})+
{1\over 3}(s_{e_3-e_5}\omega_6).\cr}$$
This shows that both ${\omega_5\over 2}$ and ${\omega_4\over 3}$ are 
non-trivial convex combinations of elements in $\oline{\Omega}$, 
hence are not extremal. It remains to show that ${\omega_2\over 2}$
is not extremal. From the definition of $\Sigma$ it is easy to check that 
${\omega_2\over 2}+t(e_1-e_2)\in \oline \Omega$ for $t\in\R$ and $|t|$ small. 
Thus ${\omega_2\over 2}$ is not extremal, concluding the 
proof of the lemma.\qed

\msk\nin {\bf The system $E_7$.} Let 
$$V\:=\{ v\in \R^8\: \la v, e_7+e_8\ra=0\}.$$
The root system $E_7$ is defined by 

$$\Sigma\:=\{ \pm e_i\pm e_j\: 1\leq i<j\leq 6\}\amalg\{\pm(e_7-e_8)\}
\amalg\{ {1\over 2}\sum_{j=1}^8\: 
(-1)^{n(j)} e_j\in V\: \sum_{j=1}^8 n(j) \ \hbox{even}\}.$$

Further a positive system is given by 
$$\eqalign{\Sigma^+\:=&\{ e_i\pm e_j\: 1\leq j<i\leq 6\}\amalg\{e_8-e_7\}\cr 
&\amalg \{ {1\over 2}(e_8-e_7 +\sum_{j=1}^6 (-1)^{n(j)} e_j)\: \sum_{j=1}^6 
n(j) \ \hbox{odd}\}.\cr}$$

Note that $\la \alpha,\alpha\ra =2$ for all $\alpha\in\Sigma$, i.e., 
all roots have the same length. 
A basis for $\Sigma^+$ is given by 

$$\eqalign{&\Pi=\{\alpha_1, \ldots,\alpha_6,\alpha_7\}\cr 
&=\{ {1\over 2}(e_8-e_7-e_6-e_5-e_4-e_3-
e_2+e_1), e_1 +e_2, e_2-e_1, e_3-e_2, e_4-e_3, e_5-e_4, e_6-e_5\}.\cr}$$
The fundamental weights 
are:

$$\eqalign{\omega_1 &= (e_8-e_7)\cr 
\omega_2 &= {1\over 2} (e_1+e_2+e_3+e_4+e_5+e_6)+ (e_8-e_7)\cr 
\omega_3 &= {1\over 2} (-e_1 + e_2 +e_3 +e_4 +e_5+e_6) +{3\over 2}(e_8-e_7)\cr 
\omega_4 &=  (e_3 +e_4 +e_5+e_6) + 2(e_8-e_7)\cr 
\omega_5 &= (e_4+e_5+e_6) +{3\over 2} (e_8-e_7)\cr 
\omega_6 &=  (e_5 +e_6)+ (e_8-e_7)\cr
\omega_7 &= e_6 +{1\over 2}(e_8-e_7).\cr}$$

The highest root $\beta$ in $\Sigma^+$ has the expression

$$\beta=2\alpha_1+2\alpha_2+3\alpha_3+ 4\alpha_4+3\alpha_5+2\alpha_6+\alpha_7.$$

\Lemma 3.18. For a root system of type $E_7$ we have 
$$\partial_e\Omega=\Ext(\oline \Omega)= {\cal W} ({\pi\over 2}\omega_7).$$

\Proof. In view of Lemma 3.16(iii), it is sufficient to 
show that ${\omega_1\over2}$, ${\omega_2\over 2}$, ${\omega_3\over 2}$, ${\omega_4\over 4}$, 
${\omega_5\over 3}$ and ${\omega_6\over 2}$  
are not extremal in $\oline {\Omega}$. Simple verification shows that 

$$\eqalign{{\omega_1\over 2}& ={1\over 2}(\omega_7) + {1\over 2}(
s_{e_5-e_6}s_{e_5+e_6} \omega_7)\cr 
{\omega_3\over 3}& ={2\over 3}({\omega_2\over 2}) + {1\over 3}(
s_{e_5-e_6}s_{e_5+e_6} \omega_7)\cr 
{\omega_4\over 4}& ={3\over 4}({\omega_5\over 3}) + {1\over 4}(
s_{e_3-e_6} \omega_7)\cr 
{\omega_5\over 3}& ={1\over 3}(\omega_7) + {2\over 3}(
{\omega_6\over 2})\cr 
{\omega_6\over 2}& ={1\over 2}(\omega_7)
 + {1\over 2}(s_{e_5-e_6}\omega_7).\cr}$$
Finally, it is easy to show that ${\omega_3\over 2}+t(e_1-e_2)\in 
\oline{\Omega}$ for $t\in\R$ small. \qed

\msk Having this information we can now handle all the exceptional 
non-compactly causal  Lie algebras. We write $E_6$ and $E_7$ 
for simply connected complex Lie groups with complex exceptional 
Lie algebra $\e_6$ resp. $\e_7$. By $E_{6(*)}$ and $E_{7(*)}$ we denote the 
real forms of $E_6$ and $E_7$ with real forms $\e_{6(*)}$, resp. $\e_{7(*)}$. 
\msk 

\msk\nin {\bf The case of $E_{6(6)}$.}
Let$\g=\e_{6(6)}$. Note that $\g$ is a normal real form 
of $\g_\C$, i.e.,  $\a$ is a Cartan algebra of $\g$. In particular 
$\Sigma$ is of type $E_6$.  

\par {} From Lemma 3.16 we obtain that 
$$\partial_e\Omega={\cal W}({\pi\over 2}\omega_1)\amalg{\cal W}(
{\pi\over 2}\omega_6).$$
From the structure of $\Sigma$ (in particular from the formula for 
the highest root) we obtain that 
$$\Spec (\ad {\pi\over 2}\omega_1)=\Spec (\ad {\pi\over 2}\omega_6)=
\{-{\pi\over 2}, 0, {\pi\over 2}\}.$$
Hence the prescriptions 
$$\sigma_1(X)=\Ad(\exp(i\pi \omega_1)(X)\qquad\hbox {and}\qquad  \sigma_2(X)=
\Ad(\exp(i \pi \omega_6))(X)$$
define two involutions on $\g$. 
Then $\tau_1\:=\sigma_1\circ \theta$ and $\tau_2\:=\sigma_2\circ \theta$
are involutions on $\g$. Note that both involutions are conjugate 
under the outer automorphism $\kappa$.
\par Recall from the list in Remark 3.3(a) that there is up 
to conjugation a unique involution $\tau$ on $\g$ 
which turns $(\g,\tau)$ into a non-compactly causal symmetric
Lie algebra. Thus Lemma 3.4 and   
Theorem 3.5 imply that:

\Proposition 3.19. If $\g=\e_{6(6)}$, then 
$$\partial_d\Xi\simeq E_{6(6)}/ E_{6(6)}^{\tau_1}\amalg E_{6(6)}/
E_{6(6)}^{\tau_2}$$
with $\g^{\tau_1}\simeq \g^{\tau_2}\simeq \sp(2,2)$.\qed

\msk\nin{\bf The case of $E_6$}.  
This case is completely analogous to the $e_{6(6)}$-case. 
One obtains that: 

\Proposition 3.20. If $\g=\e_6$, then 

$$\partial_d\Xi\simeq E_6/ E_{6(-14)}\dot \amalg E_6/E_{6(-14)}\qeddis

\msk\nin {\bf The case of $E_{6(-26)}$.} Let $\g=\e_{6(-26)}$. 
In this case $\Sigma$ is of type $A_2$ and we can use our results 
from our discussions from the root systems $A_n$. In particular we have 

$$\partial_e\Omega={\cal W}(Y_1) \amalg 
{\cal W}(Y_2).$$
Note that $Y_1$ and $Y_2$ are conjugate under the outer 
isomorphism of the Dynkin diagram. Hence we obtain 
two conjugate involutions on $\g$ by 

$$\tau_j(X)=(\Ad(\exp(i2Y_j))\circ \theta)(X)$$
for $j=1,2$. From Lemma 3.4 and Theorem 3.5 we thus 
obtain that: 

\Proposition 3.21. If $\g=\e_{6(-26)}$, then 

$$\partial_d\Xi\simeq E_{6(-26)}/ E_{6(-26)}^{\tau_1}\amalg E_{6(-26)}/
E_{6(-26)}^{\tau_2}$$
with $\g^{\tau_1}\simeq \g^{\tau_2}\simeq \f_{4(-20)}$.\qed

\msk\nin {\bf The case of $E_{7(7)}$.}
Let $\g=\e_{7(7)}$. Note that $\g$ is a normal real form 
of $\g_\C$, i.e.,  $\a$ is a Cartan algebra of $\g$. In particular 
the restricted root system is of type $E_7$. 

\par By Lemma 3.18 we have $\partial_e\Omega={\cal W}({\pi\over 2}\omega_7)$.
As before we conclude that the prescription  

$$\tau(X)=(\Ad(\exp(i\pi \omega_7))\circ\theta)(X)$$
defines an involution on $\g$. 
Hence Lemma 3.4, Theorem 3.5 and the list in Remark 3.3(a)
imply:

\Proposition 3.22. If $\g=\e_{7(7)}$, then 

$$\partial_d\Xi\simeq E_{7(7)}/ E_{7(7)}^{\tau},$$
where $\g^{\tau}\simeq \su^*(8)$.\qed

\msk\nin {\bf The case of $E_7$.} Let $\g=\e_7$. This 
case is completely analogous to the $\e_{7(7)}$-case and 
one obtains:

\Proposition 3.23. If $\g=\e_7$, then 

$$\partial_d\Xi\simeq E_7/ E_{7(-25)}.\qeddis

\subheadline{Classification of distinguished boundaries}

Taking all our results from the previous discussions together 
we have proved the following classification result:

\Theorem 3.24. Let $\g$ be a non-compactly causal Lie algebra. 

\item{(i)} If $G$ is locally  a classical group, 
then the situation is as follows:
\itemitem{(a)} If $G$ is not a special orthogonal group, 
then 
$$\partial_d\Xi\simeq \coprod_{j=1}^m G/H_j$$
where every $G/H_j$ is a non-compactly causal symmetric space. With 
$\h_j\:={\rm Lie}(H_j)$ we have: 
$$\vbox{\tabskip=0pt\offinterlineskip
\def\tablerule{\noalign{\hrule}}
\halign{\strut#&\vrule#\tabskip=1em plus2em&
\hfil#\hfil&\vrule#&\hfil#\hfil &\vrule#\tabskip=0pt\cr\tablerule
&&\omit\hidewidth $\g$\hidewidth&& 
\omit\hidewidth $\coprod_{j=1}^m\h_j$\hidewidth & \cr\tablerule
&& $\sp(n,\R)$ && $\gl(n,\R)$ & \cr\tablerule
&& $\sp(n,\C)$ && $\sp(n,\R)$ & \cr\tablerule
&& $\sp(n,n)$  && $\sp(n,\C)$ & \cr\tablerule 
&& $\su(n,n)$ && $\sL(n,\C)\oplus\R$ &  \cr\tablerule
&& $\so^*(4n)$  && $\sL(n,\H)\oplus\R$ & \cr\tablerule
&& $\sL(n,\R)$ && $\coprod_{q=1}^{n-1}\so(q,n-q)$ & \cr\tablerule
&& $\sL(n,\C)$ && $\coprod_{q=1}^{n-1}\su(q,n-q)$ & \cr\tablerule
&& $\sL(n,\H)$ && $\coprod_{q=1}^{n-1}\sp(q,n-q)$ & \cr\tablerule
}}$$

\itemitem{(b)}If $G$ is locally $\SO(n,n)$  or $\SO(2n,\C)$
for $n\geq 3$, then 
$$\partial_d\Xi\simeq G/H_1\amalg G/ H_2\amalg G/H_2$$
with $G/H_1$ and $G/H_2$ non-compactly causal symmetric spaces and: 
$$\vbox{\tabskip=0pt\offinterlineskip
\def\tablerule{\noalign{\hrule}}
\halign{\strut#&\vrule#\tabskip=1em plus2em&
\hfil#\hfil&\vrule#&\hfil#\hfil &\vrule#&
\hfil#\hfil&\vrule# \tabskip=0pt\cr\tablerule
&&\omit\hidewidth $\g$\hidewidth&& 
\omit\hidewidth $\h_1$\hidewidth&& 
\omit\hidewidth $\h_2$\hidewidth& \cr\tablerule
&& $\so(n,n)$   && $\so(n-1,1)\oplus\so(n-1,1)$ && $\so(n,\C)$ &\cr\tablerule
&& $\so(2n,\C)$ && $\so(2,2n-2)$ && $\so^*(2n)$ & \cr\tablerule}}$$
\itemitem{(c)} If $G$ is locally $\SO(p,q)$ for $1\leq p<q$  or 
$\SO(2n+1,\C)$, then for 
$p,n\geq 3$ one has  
$$\partial_d\Xi\simeq G/H_1\amalg G/ H_2$$
with $G/H_1$ non-compactly causal and $G/H_2$ a homogeneous
but non-symmetric space:
$$\vbox{\tabskip=0pt\offinterlineskip
\def\tablerule{\noalign{\hrule}}
\halign{\strut#&\vrule#\tabskip=1em plus2em&
\hfil#\hfil&\vrule#&\hfil#\hfil &\vrule#&
\hfil#\hfil&\vrule# \tabskip=0pt\cr\tablerule
&&\omit\hidewidth $\g$\hidewidth&& 
\omit\hidewidth $\h_1$\hidewidth&& 
\omit\hidewidth $\h_2$\hidewidth& \cr\tablerule
&& $\so(p,q)$ && $\so(p-1,1)\oplus \so(q-1,1)$ && $\so(p,\C)\oplus
\so(q-p)$ &\cr\tablerule
&& $\so(2n+1,\C)$ && $\so(2,2n-1)$ && $\so^*(2n)$ & \cr\tablerule}}$$
In the low-dimensional cases $p=1,2$ and $n=2$ one has 
$$\partial_d\Xi\simeq G/H$$
with $G/H$ non-compactly causal and:
$$\vbox{\tabskip=0pt\offinterlineskip
\def\tablerule{\noalign{\hrule}}
\halign{\strut#&\vrule#\tabskip=1em plus2em&
\hfil#\hfil&\vrule#&\hfil#\hfil &\vrule#\tabskip=0pt\cr\tablerule
&&\omit\hidewidth $\g$\hidewidth&& 
\omit\hidewidth $\h$\hidewidth & \cr\tablerule
&& $\so(1,q)$ && $\so(1,q-1)$ & \cr\tablerule
&& $\so(2,q)$ && $\so(1,q-1)\oplus\R $ & \cr\tablerule
&& $\so(5,\C)$ && $\so(2,3)$ & \cr\tablerule}}$$
\item{(ii)} For the exceptional cases with $\g_\C\simeq \e_6$ or 
$\g_\C \simeq \e_6\oplus\e_6$ we have $\partial_d\Xi\simeq G/H \amalg G/H$
with $G/H$ non-compactly causal. If $\h={\rm Lie}(H)$, then: 
$$\vbox{\tabskip=0pt\offinterlineskip
\def\tablerule{\noalign{\hrule}}
\halign{\strut#&\vrule#\tabskip=1em plus2em&
\hfil#\hfil&\vrule#&\hfil#\hfil &\vrule#\tabskip=0pt\cr\tablerule
&&\omit\hidewidth $\g$\hidewidth&& 
\omit\hidewidth $\h$\hidewidth & \cr\tablerule
&& $\e_{6(6)}$ && $\sp(2,2)$ & \cr\tablerule
&& $\e_{6(-26)}$ && $\f_{4(-20)}$ & \cr\tablerule
&& $\e_6$ && $\e_{6(-14)}$ & \cr\tablerule}}$$
\item{(iii)} For the exceptional cases with $\g_\C\simeq \e_7$ or 
$\g_\C \simeq \e_7\oplus\e_7$ the distinguished boundary $\partial_d\Xi\simeq G/H$
is connected. Further $G/H$ is non-compactly causal
and we have:  
$$\vbox{\tabskip=0pt\offinterlineskip
\def\tablerule{\noalign{\hrule}}
\halign{\strut#&\vrule#\tabskip=1em plus2em&
\hfil#\hfil&\vrule#&\hfil#\hfil &\vrule#\tabskip=0pt\cr\tablerule
&&\omit\hidewidth $\g$\hidewidth&& 
\omit\hidewidth $\h$\hidewidth & \cr\tablerule
&& $\e_{7(7)}$ && $\su^*(8)$ & \cr\tablerule
&& $\e_{7(-25)}$ && $\e_{6(-26)}\times\R$ & \cr\tablerule
&& $\e_7$ && $\e_{7(-25)}$ & \cr\tablerule
}}\qeddis

{}From the classification we obtain the following important 
result:

\Theorem  3.25. Let $\g$ be a non-compactly causal Lie algebra
and $\partial_d\Xi\simeq\coprod_{j=1}^n G/H_j$ the decomposition 
of $\partial_d\Xi$ into $G$-orbits.  
Then the following assertions hold: 

\item{(i)} If $G/H$ is a non-compactly causal symmetric space, 
then  $G/H$ appears locally as a component of $\partial_d\Xi$. 
\item{(ii)} A boundary component  $G/H_j$ of $\partial_d\Xi$
is totally real if and only if $G/H_j$ is symmetric. \qed 

\Remark 3.26. In the situation of Theorem 3.24
there is a boundary component of $\partial_d\Xi$
with complex dimensions if and only if 
$G$ is locally $\SO(p,q)$ for $2<p<q$ or $\SO(2n+1,\C)$ for $n>2$. 
In these cases we have $\partial_d \Xi=G/H_1\amalg G/H_2$ with 
$G/H_1$ totally real and $G/H_2$ a component admitting 
complex submanifolds of $G_\C/K_\C$.\qed

\sectionheadline{4. Relations to Jordan algebras -- the cases with 
$A_n$-type root systems}

We first recall some standard facts of Jordan algebras and explain 
some results from [KrSt01b].

\par Let $V$ be an Euclidean Jordan algebra. For $x\in V$ we define 
$L(x)\in \End(V)$ 
by $L(x)y=xy$, $y\in V$. We denote by $e$ the identity element of $V$. 
The symmetric cone $W\subeq V$ associated to $V$ can be defined as 
$$W\:=\Int\{ x^2\: x\in V\},$$
where $\Int(\cdot)$ denotes the interior of $(\cdot)$. 
We write $G\:=\Aut(W)_0$ for the connected component of the automorphism group of $W$
which contains the identity. Recall 
that $G$ is a reductive subgroup of $\Gl(V)$. The isotropy group in $e$
$$K\:=G_e\:=\{ g\in G\: g(e)=e\}$$
is a maximal compact subgroup of $G$. The mapping 
$$G/K\to W, \ \ gK\mapsto g(e)$$
is a homeomorphism. 
\par We write $V_\C =V+iV$ for the complexification  of $V$ and define the tube domain 
$$S_W\:=V+i W\subeq V_\C.$$
Let $S$ denote the connected component of the complex automorphism group of 
$S_W$
which contains the identity. Write 
$U\:=S_{ie}$ for the stabilizer of $ie\in S_W$. Then $U$ is a maximal compact subgroup 
of $S$ extending $K$ and the mapping 
$$S/U\to S_W, \ \ gU\mapsto g(ie)$$
is a biholomorphism of the Hermitian symmetric space $S/U$ onto $S_W$.
Note that $G/K\subeq S/U$.  

\par For the convenience of the reader we list here all irreducible Euclidean 
Jordan algebras and their associated groups $G$ and $S$ (cf.\ [FaKo94, p.\ 213]). 

$$\vbox{\tabskip=0pt\offinterlineskip
\def\tablerule{\noalign{\hrule}}
\halign{\strut#&\vrule#\tabskip=1em plus2em&
\hfil#\hfil&\vrule#&\hfil#\hfil &\vrule#&
\hfil#\hfil&\vrule# \tabskip=0pt\cr\tablerule
&&\omit\hidewidth $V$\hidewidth&& 
\omit\hidewidth $G$\hidewidth&& 
\omit\hidewidth $S$\hidewidth& \cr\tablerule
&& $\Symm(n,\R)$ && $\Gl(n,\R)_+$ && $\Sp(n,\R)$ &\cr\tablerule
&& $\Herm(n,\C)$ && $\Sl(n,\C)\times \R^+$ && $\SU(n,n)$ & \cr\tablerule
&& $\Herm(n,\H)$ && $\Sl(n,\H)\times \R^+$ && $\SO^*(4n)$ & \cr\tablerule
&& $\R\times\R^n$ && $\SO(1,n)\times\R^+$ && $\SO(2,n+1)$ &\cr\tablerule
&& $\Herm(3,{\Bbb O})$ && $E_{6(-26)}\times \R^+$ && $E_{7(-25)}$ & \cr\tablerule}}
$$

\par Let $c_1,\ldots, c_n$ be a Jordan frame (cf.\ [FaKo94, p.\ 44]) of $V$ and set 
$V^0\:=\bigoplus_{j=1}^n \R c_j$. Recall that $K(V^0)=V$ (cf.\ [FaKo94, Cor.\ IV.2.7])

The choice $\a\:=\bigoplus_{j=1}^n \R L(c_j)$ defines a maximal abelian 
hyperbolic subspace orthogonal to $\k$. As the table shows, the root system 
$\Sigma=\Sigma(\g,\a)$ is classical and of type $A_{n-1}$. 
If we define $\eps_j\in \a^*$ by $\eps_j(L(c_i))=\delta_{ij}$, then 
we have 
$$\Sigma=\{ {1\over 2}(\eps_i-\eps_j)\: i\neq j\}$$
(cf.\ [FaKo94, Prop.\ VI.3.3]). Thus 
$$\Omega=\{ x=\sum_{j=1}^n x_j L(c_j)\: x_j\in \R, \ |x_i-x_j|<\pi\}.$$
In particular, for  
$$\Omega_0\:=\bigoplus_{j=1}^n ]-{\pi\over 2}, {\pi\over 2}[ L(c_j)$$
we have $\Omega_0\subeq \Omega$. 
Note that the domain 
$$\Xi_0\:= G \exp(i\Omega_0)K_\C/ K_\C \subeq G_\C/ K_\C$$
is a $G$-invariant open subdomain of $\Xi$ (cf.\ [KrSt01, Lemma 1.4]).

\Theorem 4.1. The mapping 
$$\Xi_0\to S_W,  \ \ gK_\C\mapsto g(ie)$$
is a biholomorphism. 

\Proof. [KrSt01b, Th.\ 2.2.].\qed 

\subheadline{The distinguished boundary of $\Xi_0$}

As in Section 2 we prove the following result:

\Proposition 4.2. The distinguished boundary of $\Xi_0$ in $G_\C/ K_\C$ is 
given by

$$\partial_d\Xi_0=G\exp(i\partial_e\a_1)K_\C/ K_\C$$
with

$$\partial_e\Omega_0=\{ {\pi\over 2}(\pm L(c_1)\pm\ldots\pm L(c_n))\}.\qeddis

Recall that the Weyl group ${\cal W}$ is isomorphic to $S_n$ and 
acts as the full permutation group of $L(c_1), \ldots, L(c_n)$. 
Thus we immediately obtain that: 

\Lemma 4.3. We have the disjoint union 

$$\partial_e\Omega_0=\coprod_{p=0}^n {\cal W} (Y_p)$$
with 
$$Y_p\:={\pi\over 2}(L(c_1)+\ldots + L(c_p)- L(c_{p+1}) - \ldots - L(c_n)).
\qeddis

Recall that $G_\C/K_\C$ sits canonically in $V_\C$ via 
the embedding 

$$G_\C/K_\C\to V_\C, \ \ gK_\C\mapsto g(ie).$$
In the sequel we identify $G_\C/ K_\C$ as a subset of $V_\C$. 

Set $z_p\:=\exp(iY_p)K_\C\in V_\C$. Then we have:

\Lemma 4.4. For all $0\leq p\leq n$ we have that 

$$z_p=c_1+\ldots+ c_p - c_{p+1}- \ldots -c_n.\qeddis

Write $G_p$ for the stabilizer of $z_p$ in $G$.

\Example 4.5. We take $V=\Symm(n,\R)$ with $G=\GL(n,\R)_+$. The action 
of $G$ is given by 

$$G\times V\to V, \ \ (g,X)\mapsto gXg^t.$$
Further we have that

$$z_p=\diag (1,\ldots 1, -1,\ldots, -1)$$
with $p$-times $+1$ on the diagonal. Thus we see that 

$$(\forall 0\leq p\leq n) \qquad G_p=\SO(p,n-p).$$
In particular $G_0=G_n=\SO(n,\R)$ is compact and all other 
$G_p$ are non-compact. 
Write $V_p$ in $V$ for the symmetric matrices with signature $(p, n-p)$.
Further we write $V^{\rm reg}$ for the invertible matrices in $V$
and set $V_p^{\rm reg}\:=V_p\cap V^{\rm reg}$. Then we have 

$$V^{\rm reg}=\coprod_{p=0}^n V_p^{\rm reg}$$
and 
$$  V_p^{\rm reg} =G(z_p)\simeq G/G_p.$$
In particular we see that $\partial_d\Xi_0$ is a Zariski open 
subset of $V$. \qed 

We are now going to generalize the results in Example 4.5 to 
arbitrary irreducible Euclidean Jordan algebras. 

Let $V$ be an irreducible Euclidean Jordan algebra. 
Write $\det(x)$ for the Jordan algebra determinant of $V$ and note that 
$\det$ is a polynomial function on $V$. Define the subset of 
regular elements of $V$ by 

$$V^{\rm reg}\:=\{ x\in V\: \det x\neq 0\}$$
and note that $V^{\rm reg}$ is a Zariski open subset of $V$. 

For every $x\in V$ we define a real polynomial 

$$f(\lambda, x)\:=\det(\lambda e -x).$$
Note that for fixed $x$, the polynomial $f(\lambda, x)$ is completely 
reducible over $\R$. If $x$ is regular, then it has degree $n$ and 
all roots are non-zero. 
If $f(\lambda)$ is a polynomial, then we define its 
signature $\sgn f$ to be the number of its positive roots. 
For every $0\leq p\leq n$ now set 

$$V_p\:=\{ x\in V\: \sgn f(\cdot, x) =p\}.$$ 

Note that 

$$V^{\rm reg}=\coprod_{p=0}^n V_p^{\rm reg}$$
is a disjoint decomposition in cones. The only convex cones 
in this decomposition are $V_0^{\rm reg}=-W$ and $V_n^{\rm reg}=
W$. 
\par Observe that $z_p\in V_p$ for all $0\leq p\leq n$.

\Proposition 4.6. For every $0\leq p\leq n$ the set $V_p^{\rm reg}$ is 
$G$-invariant and $G$ acts transitively on it.  Hence the map

$$G/G_p\to V_p^{\rm reg}, \ \ gG_p\mapsto g(z_p)$$
is an isomorphism. 

\Proof. First we show that $V_p^{\rm reg}$ is invariant under $K$. 
For that fix $x\in V_p^{\rm reg}$ and $k\in K$. Recall that $k(e)=e$
and that $K$ acts on $V$ by Jordan algebra 
automorphisms. We have 

$$f(\lambda, k(x))=\det (\lambda e-k(x))=\det(k(\lambda(e)-x))
=\det(\lambda e-x)$$
where in the last equality we used the fact that $\det k(y)=\det y$
for all $y\in V$ and $k\in K$.  
\par It is straightforward to check that 

$$V_p^{\rm reg}\cap V^0 =A(z_p).$$
Since $V_p^{\rm reg}$ is $K$-invariant and $G=KAK$ the 
assertions of the proposition now follow.\qed 

\Corollary 4.7. The distinguished boundary $\partial_d\Xi_0$ 
,realized in $V_\C$, is a Zariski open subset of $V$.\qed 

Next we compute the isotropy groups $G_p$ for $0\leq p\leq n$. 

\Proposition 4.8. For an irreducible Euclidean Jordan algebra 
$V$ the isotropy groups $G_p$ for $0\leq p\leq n$ are given as 
follows: 

\item{(i)}  For the classical matrix Jordan algebras 
one has: 

$$\vbox{\tabskip=0pt\offinterlineskip
\def\tablerule{\noalign{\hrule}}
\halign{\strut#&\vrule#\tabskip=1em plus2em&
\hfil#\hfil&\vrule#&\hfil#\hfil &\vrule#&
\hfil#\hfil&\vrule# \tabskip=0pt\cr\tablerule
&&\omit\hidewidth $V$\hidewidth&& 
\omit\hidewidth $G$\hidewidth&& 
\omit\hidewidth $G_p$\hidewidth& \cr\tablerule
&& $\Symm(n,\R)$ && $\Gl(n,\R)_+$ && $\SO(p,n-p)$ &\cr\tablerule
&& $\Herm(n,\C)$ && $\Sl(n,\C)\times \R^+$ && $\SU(p,n-p)$ & \cr\tablerule
&& $\Herm(n,\H)$ && $\Sl(n,\H)\times \R^+$ && $\Sp(p,n-p)$ & \cr\tablerule}}
$$
\item{(ii)} For $V=\R\times \R^n$ and $G=\SO(1,n)\times\R^+$ one has: 
\itemitem{(a)} $G_p=\SO(n)$ for $p=0,2$. 
\itemitem{(b)} $G_p=\SO(1,n-1)$ for $p=1$. 
\item{(iii)} For $V=\Herm(3,{\Bbb O})$ and $E_{6(-26)}\times\R^+$ one has: 
\itemitem{(a)} $G_p=K=F_{4(-52)}$ for $p=0,3$. 
\itemitem{(b)} $G_p\simeq F_{4(-20)}$ for $p=1,2$.

\Proof. (i) For all classical matrix algebras this is the same 
computation as in Example 4.5.
\par\nin (ii) Straightforward calculation similar 
to the one in (i).
\par\nin (iii) While (a) is clear, (b) is shown as in 
Proposition 3.21.
\qed

\Corollary 4.9. The distinguished boundary $\partial_d\Xi_0$ is $G$-isomorphic 
to $\partial_d\Xi_0=\coprod_{p=0}^n G/G_p$ where

\item{(i)} $G/G_p$ is a Riemannian symmetric space for $p=0, n$. 
\item{(ii)} $G/G_p$ is a non-compactly causal symmetric space 
for $p\neq 0, n$.\qed 

\Remark 4.10. (a) Comparing Corollary 4.9 to our results for the 
distinguished boundary of the bigger domain $\Xi$ we see
that the distinguished boundary of 
$\Xi_0$ is the distinguished boundary of $\Xi$ plus 
two copies of the Riemannian symmetric space. Also observe that 
every boundary component of $\partial_d\Xi_0$ is totally real. 

\par\nin (b) The observation in (a) fits into the following 
general philosphy: Non-compactly causal symmetric Lie algebras $(\g,\tau)$ 
are the class of symmetric Lie algebras  which are closest to 
non-compactly Riemannian 
symmetric algebras $(\g,\theta)$. In fact these two classes make up 
the class of symmetric Lie algebras $(\g,\tau)$ where $\q$ admits 
a non-trivial open $\Ad(H)$-invariant hyperbolic convex set. 
These Lie algebras were subject of systematic study in [KrNe96].
\qed

\sectionheadline{5. Further results on the special orthogonal groups}

The examples discussed in Section 4 fit into the broader context 
of comparing $\Xi$ with symmetric spaces of Hermitian type 
$S/U$ which contain $G/K$ as a totally real submanifold 
(cf.\ [KrSt01b]). 
\par In this section we will focus on the special orthogonal 
groups $G=\SO_e(p,q)$ and $G=\SO(n,\C)$ which in some sense are 
the classical groups with the most complicated structure of $\Xi$
(cf.\ Theorem 3.24).
\par Our choice of the maximal compact subgroup $K<G$ is as in 
Section 3. If $G=\SO_e(p,q)$, then we take $S=\SU(p,q)$ and 
$U=S(U(p)\times U(q))$ and if $G=\SO(n,\C)$, then we 
choose $S=\SO^*(2n)$ and $U=U(n)$.  
The embedding 

$$G/K\into S/U, \ \ gK\mapsto gU\leqno(5.1)$$
realizes $G/K$ as a totally real submanifold of the 
symmetric space of Hermitian type $S/U$. 
\par Let $\a\subeq \p$ be a maximal abelian subspace
as in Section 3. Write $\hat\Sigma\:=\Sigma(\a,\s)$ for the 
doubly restricted root system of $\s$ with respect to $\a$. 
Notice that $\hat\Sigma$ is of type $C_n$ if 
$\g=\so(n,n)$ or $\g=\so(2n,\C)$ and of type 
$BC_n$ otherwise. In the notation of Section 3 we have:

$$\hat\Sigma=\cases{ \{ \pm\eps_i\pm\eps_j\: 1\leq i, j\leq n\}\bs\{0\}
& for $\g=\so(n,n), \so(2n,\C)$,\cr 
\{ \pm\eps_i\pm\eps_j\: 1\leq i, j\leq p\}\bs \{0\} \amalg\{ 
\pm \eps_i\: 1\leq i\leq p\} & for $\g=\so(p,q), \so(2p+1,\C)$
, $(p< q)$.\cr}$$
Define $\Omega_0\:=\{X\in\a\: (\forall \alpha\in \hat\Sigma)\ 
|\alpha(X)|<{\pi\over 2}\}$. Note that $\Omega_0\subeq \Omega$
since $\hat\Sigma\supeq \Sigma$. Define 
$$\Xi_0\:=G\exp(i\Omega_0)K_\C/K_\C$$
and notice that $\Xi_0$ is a $G$-invariant open subdomain
of $\Xi$. Finally define the distinguished boundary 
of $\Xi_0$ by $\partial_d\Xi_0\:=G\exp(i\partial_e\Omega_0)K_\C/K_\C$.

\par Write $\uu$ for the Lie algebra of $U$ and 
$\s_\C =\p^+\oplus\uu_\C\oplus\p^-$ for the Harish-Chandra
decomposition of $\s_\C$. Let ${\cal D}\subeq \p^+$ be the 
Harish-Chandra realization of $S/U$ as a bounded symmetric domain. 
The embedding in (5.1) extends to an embedding 
$$ G_\C/ K_\C\into S_\C /U_\C P^-\leqno(5.2)$$
and it follows from [KrSt01b, Th.\ 2.5] that the image 
of $\Xi_0$ under (5.2) is precisely $S/U\simeq {\cal D}$, i.e., 
we have a $G$-equivariant biholomorphism
$$\Xi_0\simeq S/U.$$

\Theorem 5.1. Let $G=\SO_e(p,q)$ or $G=\SO(n,\C)$. Then 
the distinguished boundary $\partial_d\Xi_0$ is given as follows:
\item{(i)} If $G=\SO_e(n,n)$ or $G=\SO(2n,\C)$, then 
$\partial_d\Xi_0=G/H\amalg G/H$ with 
$\h\:={\rm Lie}(H)$ given by 
$$\h= \cases{ \so(n,\C)  & for $\g=\so(n,n)$\cr 
\so^*(2n) & for $\g=\so(2n,\C)$.\cr}$$
\item{(ii)} If $G=\SO_e(p,q)$ for $p<q$ or $G=\SO(2n+1,\C)$, then 
$\partial_d\Xi_0=G/H$ with 
$\h\:={\rm Lie}(H)$ given by 
$$\h= \cases{ \so(p,\C)\oplus \so(q-p)  & for $\g=\so(p,q)$\cr 
\so^*(2n) & for $\g=\so(2n+1,\C)$.\cr}$$

\Proof. (i) Keep the notation of Lemma 3.11. Then it 
follows from the structure of $\hat\Sigma$ that 
$$\partial_e\Omega_0={\cal W}(Y_2)\amalg {\cal W}(Y_3).$$
Now the assertion follows from the computations 
before Proposition 3.12 and Proposition 3.14. 
\par\nin (ii) This is analogous to (i). 
\qed

\def\entries{

\[AkGi90 Akhiezer, D.\ N., and S.\ G.\ Gindikin, {\it On Stein extensions of 
real symmetric spaces}, 
Math.\ Ann.\ {\bf 286}, 1--12, 1990

\[BHH01 Burns, D., S.\ Halverscheid, and R.\ Hind, {\it The Geometry 
of Grauert Tubes and Complexification of Symmetric Spaces}, 
preprint

\[FaKo94 Faraut, J., and A. Koranyi, ``Analysis on symmetric cones", 
Oxford Mathematical Monographs, Oxford University Press, 1994

\[Gi98 Gindikin, S., {\it Tube domains in Stein symmetric spaces}, 
Positivity in Lie theory: open problems, 81--97, de Gruyter Exp. Math., 
{\bf 26}, de Gruyter, Berlin, 1998

\[GiKr01 Gindikin, S., and B.\ Kr\"otz, {\it Invariant Stein domains in 
Stein symmetric spaces and a non-linear complex 
convexity theorem}, MSRI preprint 2001-032

\[GK\'O01  Gindikin, S.,B.\ Kr\"otz, and G.\ \'Olafsson, {\it Hardy 
spaces for non-compactly causal symmetric spaces and the most 
continuous spectrum I: causally symmetric triples}, 
MSRI preprint 

\[GiMa01 Gindikin, S., and T. Matsuki, {\it Stein Extensions
of Riemann Symmetric Spaces and Dualities of Orbits on Flag 
Manifolds}, MSRI preprint 2001-028

\[Hi\'Ol96 Hilgert, J.\ and 
G.\ \'Olafsson, ``Causal Symmetric Spaces, Geometry and
Harmonic Analysis,'' Acad. Press, 1996 

\[Kn96 Knapp, A., ``Lie Groups Beyond an Introduction'', Birkh\"auser, 
Basel, 1996

\[KrNe96 Kr\"otz, B., and K.-H.\ Neeb, {\it On Hyperbolic Cones 
and Mixed Symmetric Spaces}, J. Lie Theory {\bf 6}, 69--146, 1996 

\[KrSt01a Kr\"otz, B., and R.J. Stanton, {\it Holomorphic extensions of representations: (I)  
automorphic functions}, Ohio-State preprint

\[KrSt01b Kr\"otz, B., and R.J. Stanton, {\it Holomorphic extensions of representations: (II)  
geometry and harmonic analysis}, preprint
}

{\sectionheadline{\bf References}
\frenchspacing
\entries\par}
\dlastpage 
\bye
\end